\documentclass[10pt]{article}
\usepackage[utf8]{inputenc}
\usepackage[margin = 1.8 cm]{geometry}

\usepackage[shortlabels]{enumitem}
\usepackage{amsthm}
\usepackage{amsmath}
\usepackage{amssymb}
\usepackage{setspace}
\usepackage{mathtools}
\usepackage{verbatim}
\usepackage{csquotes}
\usepackage{graphicx}
\usepackage[hidelinks]{hyperref}
\usepackage{thm-restate}
\usepackage{cleveref}
\usepackage{graphicx}
\usepackage{mathrsfs}
\usepackage{enumitem}
\usepackage{framed}
\usepackage{subcaption}
\usepackage{crossreftools}
\usepackage{floatrow}
\usepackage[T1]{fontenc}
\usepackage{appendix}

\theoremstyle{plain}

\newtheorem*{thm*}{Theorem}
\newtheorem{theorem}{Theorem}
\Crefname{theorem}{Theorem}{Theorems}

\newtheorem*{lem*}{Lemma}
\newtheorem{lemma}[theorem]{Lemma}
\Crefname{lemma}{Lemma}{Lemmas}

\newtheorem*{claim*}{Claim}
\newtheorem{claim}{Claim}
\crefname{claim}{Claim}{Claims}
\Crefname{claim}{Claim}{Claims}

\newtheorem{prop}[theorem]{Proposition}
\Crefname{prop}{Proposition}{Propositions}

\crefname{corollary}{Corollary}{Corollaries}

\crefname{conj}{Conjecture}{Conjectures}

\newtheorem*{conj*}{Conjecture}

\Crefname{qn}{Question}{Questions}
\newtheorem*{qn*}{Question}

\Crefname{obs}{Observation}{Observations}

\Crefname{ex}{Example}{Examples}

\theoremstyle{definition}
\newtheorem{prob}[theorem]{Problem}
\Crefname{prob}{Problem}{Problems}

\newtheorem{defn}[theorem]{Definition}
\Crefname{defn}{Definition}{Definitions}

\newtheorem*{defn*}{Definition}

\theoremstyle{remark}

\renewenvironment{proof}[1][]{\begin{trivlist}
\item[\hspace{\labelsep}{\bf\noindent Proof#1.\/}] }{\qed\end{trivlist}}

\newcommand{\ceil}[1]{
    \left\lceil #1 \right\rceil
}
\newcommand{\floor}[1]{
    \left\lfloor #1 \right\rfloor
}

\newcommand{\eps}{\varepsilon}

\renewcommand{\P}{\mathbb{P}}
\newcommand{\E}{\mathbb{E}}
\newcommand{\N}{\mathbb{N}}
\newcommand{\R}{\mathbb{R}}

\DeclareMathOperator{\bin}{Bin}
\DeclareMathOperator{\Abs}{Abs}
\DeclareMathOperator{\polylog}{polylog}

\expandafter\def\expandafter\normalsize\expandafter{%
    \normalsize
    \setlength\abovedisplayskip{8pt}
    \setlength\belowdisplayskip{8pt}
    \setlength\abovedisplayshortskip{4pt}
    \setlength\belowdisplayshortskip{4pt}
}

\usepackage[square,sort,comma,numbers]{natbib}
 \setlist[itemize]{leftmargin=*}

\makeatletter
\newcommand{\optionaldesc}[2]{%
  \phantomsection
  #1\protected@edef\@currentlabel{#1}\label{#2}%
}
\makeatother

\title{Edge-disjoint cycles with the same vertex set}
\author{Debsoumya Chakraborti\thanks{
Mathematics Institute, University of Warwick, Coventry, United Kingdom. Research supported by the European Research Council (ERC) under the European Union Horizon 2020 research and innovation programme (grant agreement No.\ 947978). 
E-mail: {\tt
\{debsoumya.chakraborti, richard.montgomery\}@warwick.ac.uk}.} 
\and Oliver Janzer\thanks{Department of Pure Mathematics and Mathematical Statistics, University of Cambridge, United Kingdom. Research supported by a fellowship at Trinity College. Email: {\tt oj224@cam.ac.uk}}
\and Abhishek Methuku\thanks{Department of Mathematics, ETH Z\"urich,
Switzerland. Research supported by SNSF grant 200021\_196965.
E-mail: {\tt
abhishekmethuku@gmail.com}} 
\and Richard Montgomery\footnotemark[1]
}
\date{}

\begin{document}

\maketitle

\begin{abstract}
In 1975, Erd\H{o}s asked for the maximum number of edges that an $n$-vertex graph can have if it does not contain two edge-disjoint cycles on the same vertex set.
It is known that Tur\'an-type results can be used to prove an upper bound of $n^{3/2+o(1)}$. However, this approach cannot give an upper bound better than $\Omega(n^{3/2})$.
We show that, for any $k\geq 2$, every $n$-vertex graph with at least $n \cdot \polylog(n)$ edges contains $k$ pairwise edge-disjoint cycles with the same vertex set, resolving this old problem  in a strong form up to a polylogarithmic factor.
The well-known construction of Pyber, R\"odl and Szemer\'edi of graphs without $4$-regular subgraphs shows that there are $n$-vertex graphs with $\Omega(n\log \log n)$ edges which do not contain two cycles with the same vertex set, so the polylogarithmic term in our result cannot be completely removed.

Our proof combines a variety of techniques including sublinear expanders, absorption and a novel  tool for regularisation, which is of independent interest. Among other applications, this tool can be used to regularise an expander while still preserving certain key expansion properties.
\end{abstract}

\section{Introduction}
Many questions in Extremal Graph Theory concern finding cycles, for example, see~\cite{corradi1963maximal, egawa1996vertex, verstraete2000arithmetic}, and the nice survey by Verstra\"ete~\cite{verstraete2016extremal}. There have been numerous powerful methods for embedding cycles (and paths) developed in the past three decades, such as Robertson and Seymour’s work on graph linkage~\cite{robertson1995graph}, Krivelevich and Sudakov’s use of Depth First Search~\cite{krivelevich2013phase}, and the use of expanders in a long line of work by Krivelevich (see, e.g.,  the survey~\cite{krivelevich2019expanders}).

Determining, for each $k\geq 2$, the maximum number of edges an $n$-vertex graph can have without containing a cycle with length $2k$ is a fundamental problem in Tur\'an theory that goes back to a paper by Erd\H{o}s in 1938~\cite{erdos1938sequences}. The famous conjecture that the answer to this problem should be $\Omega(n^{1+\frac{1}{k}})$, to match an upper bound of Bondy and Simonovits~\cite{bondy1974cycles}, remains wide open for $k = 4$ and $k \ge 6$.
What can we say about the cycles in graphs with $n^{1+o(1)}$ edges? This has long been a subject of interest, and, for example, recently, answering a question of Erd\H{o}s~\cite{erdos1984some} from 1984, Liu and Montgomery~\cite{liu2023solution} showed among other related results that there is a constant $c$ such that any $n$-vertex graph with at least $cn$ edges contains a cycle whose length is a power of $2$.

In 1975, Erd\H{o}s~\cite[Problem 29]{erdHos1975problems} proposed the following three simple-looking problems concerning cycles. What is the maximum number of edges that an $n$-vertex graph can have if it does not contain

\hspace{1cm}\begin{minipage}{0.75\textwidth}
\begin{itemize}[noitemsep,nolistsep]
    \item  two edge-disjoint nested cycles?
    \item  two edge-disjoint nested cycles with no geometric crossings?
    \item  two edge-disjoint cycles with the same vertex set?
\end{itemize}
\end{minipage}

Here, a sequence of cycles $C_1,C_2,\ldots,C_k$ is said to be \emph{nested} if $V(C_1)\subseteq V(C_2)\subseteq\ldots \subseteq V(C_k)$.
Two edge-disjoint nested cycles are said to have
\emph{no geometric crossings} if, when the longer cycle is embedded to make a convex polygon, the other cycle has no crossing edges.

Bollob\'as~\cite{Bollobas} resolved the first problem in 1978 by showing that there is a constant $c$ such that any $n$-vertex graph with at least $cn$ edges contains two edge-disjoint nested cycles. In 1996, Chen, Erd\H{o}s and Staton~\cite{chen1996proof} extended this result by showing there is a constant $c_k$ such that any $n$-vertex graph with at least $c_k n$ edges contains $k$ edge-disjoint nested cycles.
Fern{\'a}ndez, Kim, Kim and Liu~\cite{fernandez2022nested} resolved the second problem by showing that there is a constant $c$ such that any $n$-vertex graph with at least $cn$ edges must contain two edge-disjoint nested cycles with no geometric crossings. 

The third problem, however, is different because the answer is not linear in $n$. Indeed, Pyber, R{\"o}dl, and Szemer{\'e}di~\cite{pyber1995dense} observed that this problem is related to the well-known Erd\H{o}s-Sauer problem (as discussed further later), and used their remarkable construction of $n$-vertex graphs with $\Omega(n \log\log n)$ edges and no $4$-regular subgraph, to show that there are $n$-vertex graphs with $\Omega(n \log\log n)$ edges which do not contain two edge-disjoint cycles with the same vertex set. In this paper, we consider this third problem of Erd\H os, which we restate for convenience.

\begin{prob}
[Erd\H{o}s~\cite{erdHos1975problems}, 1975]
\label{prob:erdosprob1975}
What is the maximum number of edges that an $n$-vertex graph can have if it does not contain two edge-disjoint cycles with the same vertex set?
\end{prob}

In the nearly 50 years since Erd\H{o}s posed this problem, it has been reiterated several times in the literature, including by Bollob\'{a}s \cite[Chapter~7, Problem~23]{bollobás1978extremal} in his popular book \emph{Extremal Graph Theory}, by Pyber, R{\"o}dl and Szemer{\'e}di~\cite{pyber1995dense} in 1995 in the context of the Erd\H{o}s-Sauer problem as mentioned above, and by
Chen, Erd{\H{o}}s and Staton~\cite{chen1996proof} in 1996 in the context of finding cycles with many chords.
However, despite this, and the development of powerful embedding techniques for paths and cycles, the only progress on Problem~\ref{prob:erdosprob1975} has been concluded from results on Tur\'an numbers. Here,
Chen, Erd{\H{o}}s and Staton~\cite{chen1996proof} observed that the upper bound $O(n^{7/4})$ for Problem~\ref{prob:erdosprob1975} follows from a well-known theorem of K\H{o}v\'{a}ri, S\'{o}s and Tur\'{a}n~\cite{kHovari1954problem}, since the complete bipartite graph $K_{4,4}$ contains two edge-disjoint cycles with the same vertex set.
More generally, we can use as an upper bound the Tur\'an number of any fixed graph $H$ which contains two edge-disjoint cycles with the same vertex set. This approach, however, cannot yield an upper bound better than $n^{3/2+o(1)}$, by a simple application of the probabilistic deletion method. In recent years it has been observed~\cite{Personal22} that this upper bound can indeed be proven using developments on Tur\'an-type problems. For instance, as the $2$-blowup of a cycle contains two edge-disjoint cycles with the same vertex set, a result of Janzer~\cite{janzer2023rainbow} gives an upper bound of $n^{3/2+o(1)}$ for Problem~\ref{prob:erdosprob1975}. (Here the $2$-blowup of a graph $F$ is the graph obtained by replacing each vertex of $F$ with an independent set of size $2$ and each edge of $F$ by a $K_{2,2}$.)

In this paper, we give the first methods to directly approach this problem (Problem~\ref{prob:erdosprob1975}), greatly reducing the gap between the previous best known upper bound of $n^{3/2+o(1)}$ and the best known lower bound of $\Omega(n\log\log n)$, and resolving the problem up to a polylogarithmic factor in the following stronger form.

\begin{theorem} \label{thm:main} There is some $t$ such that the following holds.
    For each $k \ge 2$, there is a constant $c = c(k)$ such that any $n$-vertex graph with at least $cn(\log n)^t$ edges contains $k$ pairwise edge-disjoint cycles with the same vertex set.
\end{theorem}

Key to our result is a novel technique for efficiently finding a nearly regular subgraph in a graph whose vertex degrees differ by at most some constant multiple.  Importantly, this technique combines well with sublinear expansion properties (see Section~\ref{sec: sublinear expander}), allowing us to find subgraphs which are very close to being regular and which have some expansion properties. Our new regularisation method has further applications, as discussed in Section~\ref{sec:concluding}, and is thus of independent interest.

\textbf{More on lower bounds and further related results.}
As mentioned above, Pyber, R{\"o}dl, and Szemer{\'e}di~\cite{pyber1995dense} noted that the problem of finding edge-disjoint cycles on the same vertex set is related to the well-known Erd\H{o}s-Sauer problem~\cite{erdos1975some} which asks for the maximum number of edges an $n$-vertex graph can have without containing a $k$-regular subgraph.  Indeed, for any $k \ge 2$, the graph consisting of $k$ pairwise edge-disjoint cycles on the same vertex set is a $2k$-regular graph.
Answering a question of Erd\H os and Sauer, and building on important work of Alon, Friedland and Kalai \cite{alon1984regular}, Pyber \cite{pyber1985regular} famously proved that for each $\eps>0$ and positive integer $k$, if $n$ is sufficiently large, then every $n$-vertex graph with at least $n^{1+\eps}$ edges contains a $k$-regular subgraph.
The construction of Pyber, R{\"o}dl, and Szemer{\'e}di~\cite{pyber1995dense} mentioned before shows that, in fact, there is a constant $c > 0$ such that for every $n$, there exists an $n$-vertex graph with at least $cn \log \log n$ edges which does not contain a $k$-regular subgraph for any $k \ge 3$. Janzer and Sudakov~\cite{janzer2023resolution} recently showed that this lower bound is best possible, thus resolving the Erd\H{o}s-Sauer problem. More precisely, they showed that there is a constant $c=c(k)$ such that any $n$-vertex graph with at least $c n \log \log n$ edges contains a $k$-regular subgraph. It would be very interesting to determine whether our bound in \Cref{thm:main} can also be improved to $O_k(n\log \log n)$.

Note that if we have two edge-disjoint cycles $C_1, C_2$ on the same vertex set, then $C_1$ is a cycle with at least as many chords as it has vertices (as the edges of $C_2$ are chords of $C_1$).
In 1996, Chen, Erd\H{o}s and Staton~\cite{chen1996proof} considered the problem of finding such a cycle. They showed that any $n$-vertex graph with at least $2n^{3/2}$ edges has a cycle with as many chords as vertices. Improving this bound, Dragani\'c, Methuku, Munh\'a Correia and Sudakov~\cite{draganic2023cycles} recently showed that any $n$-vertex graph with at least $n (\log n)^8$ edges contains a cycle with as many chords as vertices. (Note that \Cref{thm:main} implies this result up to a polylogarithmic factor.)
 Furthermore, resolving a problem of Erd\H{o}s~\cite{erdos1975some} from 1975, Brada\v{c}, Methuku and Sudakov~\cite{bradavc2023extremal}
showed that there exists a constant $c$ such that every $n$-vertex graph with at least $cn^{3/2}$ edges contains a cycle with all diagonals (where a diagonal in a cycle is a chord joining two vertices at maximum distance on the cycle).

\textbf{Organisation.} The rest of this paper is organised as follows. In \Cref{sec:overview}, we give an overview of the proof of \Cref{thm:main}. In  Section~\ref{sec:key lemmas}, we present the key lemmas that we will need, some of whose proofs are given in the later sections.  In \Cref{sec:regularizing lemma}, we prove our key regularisation lemma. In \Cref{sec:auxiliary lemmas}, we prove some more auxiliary lemmas that are needed in our proof. We prove our main result (\Cref{thm:main}) in \Cref{sec: main proof}. In \Cref{sec: sublinear expander}, building on the work of \cite{bucic2022towards}, we prove a general lemma for expanders which we use to connect pairs of vertices using vertex-disjoint paths through random subsets of vertices. Some further upcoming applications and extensions of our regularisation lemma are discussed in \Cref{sec:concluding}. Finally, in an appendix, we include the variant we need of a standard approach to find a sublinear expander in an arbitrary graph.

\textbf{Notation.} We use standard graph theoretic notation throughout the paper. In particular, for a graph $G$, we denote by $d(G)$ its average degree, and by $\delta(G)$, $\Delta(G)$ its minimum degree and maximum degree, respectively.
A graph $G$ is called \textit{$K$-almost-regular} if $\Delta(G) \le K \delta(G)$.
We call a graph \textit{$(d\pm d')$-nearly-regular} if every vertex in the graph has degree between $d-d'$ and $d+d'$.
For a set $F\subseteq E(G)$, we denote by $G - F$ the graph obtained from $G$ after deleting all the edges in $F$.
For a set $S \subseteq V(G)$, let $G[S]$ denote the subgraph of $G$ induced by $S$. We denote the number of edges of $G$ by $e(G)$, and for $S\subseteq V(G)$, let $e_G(S)$ denote the number of edges of $G$ induced by $S$. For $v\in V(G)$ and $S\subset V(G)$, we write $d_G(v,S)$ for $|N_G(v)\cap S|$. For two disjoint sets $A, B \subseteq V(G)$, $G[A,B]$ is the bipartite subgraph of $G$ with parts $A$ and $B$ consisting of all edges of $G$ between $A$ and $B$. Moreover, $e_G(A, B)$ is the number of edges of $G$ which are incident to both $A$ and $B$.  For a bipartite graph $G$ with the bipartition $A \cup B$, a set $S\subset V(G)$ is called \emph{balanced} (with respect to the bipartition of $G$) if $|S\cap A|=|S\cap B|$. A $p$-random subset of a set $S$ is obtained by keeping each element of $S$ independently at random with probability $p$.

For distinct vertices $u, v$, a $u$-$v$ path is a path joining $u$ and $v$. For a $u$-$v$ path $P$, the vertices $u, v$ are called the \emph{endpoints} of $P$, and the rest of the vertices of $P$ are called \emph{internal vertices} of $P$.
For a set $V$, a path \emph{through} $V$ is one whose internal vertices are all in $V$.

 The asymptotic notation $o(1)$ denotes a function that tends to $0$ as $n \rightarrow \infty$, unless specified otherwise. Throughout the paper we often omit the rounding functions of real numbers for the clarity of presentation. For a positive integer $n$, we write $[n]$ to denote the set $\{1,2,\dots,n\}$. Logarithms are with base~$2$.

\section{Proof sketch} \label{sec:overview}

In order to first introduce the main ideas in our proof, we will describe how to approach a simplified version of our problem in Subsection~\ref{sec:halfproblem}, before covering the additional ideas needed for the full problem in Subsection~\ref{sec:fullproblem}. We then give a detailed outline of our proof of Theorem~\ref{thm:main} in Subsection~\ref{sec:detailedoutline} and a sketch of our key near-regularisation lemma in Subsection~\ref{sec:newreg}. In our discussion, we say a graph is `almost-regular' if its vertex degrees are the same up to a constant factor, while it is `(very) nearly regular' if it is much closer to being regular and the vertex degrees are the same up to some small additive polylogarithmic term.

\subsection{A simplified problem: combining near-regularity and expansion}\label{sec:halfproblem}
In order to introduce some of the ideas we use, we will first consider a simpler `approximate' version of the problem. Suppose we have an $n$-vertex graph $G$ with at least $n(\log n)^C$ edges (for some large constant $C$) in which we wish to find $k$ edge-disjoint cycles $C^i$, $i\in [k]$, such that $|V(C^i)\triangle V(C^j)|=o(|C^i|)$ for each $i\neq j$. Suppose we found within $G$ a $d$-regular subgraph $H$ with $d\approx (\log n)^{C/2} $ (this can be achieved with known techniques, using, e.g., Pyber~\cite{pyber1985regular}). Taking $t=(\log n)^{C/4}$, we could divide $V(H)$ into $t$ sets $V_1, V_2, \ldots,V_t$ of roughly equal size and conclude using a simple Chernoff bound that, with high probability, every $H[V_j,V_{j+1}]$ is approximately regular. Within each $H[V_j,V_{j+1}]$ we could use this near-regularity (via, for example, Vizing's theorem) to find $k$ large edge-disjoint matchings, $M_j^i$, $i\in [k]$. These matchings will be large enough that, for each $i\in [k]$, letting $F^i$ be $\cup_{j}M_j^i$ together with some isolated vertices, $F^i$ will be a linear forest with vertex set $V(H)$ with $(1+o(1))|V(H)|/t$ paths. Dividing into random vertex sets and finding matchings to glue together to create linear forests is a well known technique, used, for example in  \cite{kelmans2001asymptotically}. If we can then, edge-disjointly, join up each linear forest $F^i$ into a cycle $C^i$ using short connecting paths, we will have found $k$ edge-disjoint cycles with almost the same vertex set.

A common way to join vertices by paths is to use expansion properties. We will use a very weak form of expansion called \emph{sublinear expansion} (which originates from the work of Koml\'os and Szemer\'edi \cite{komlos1994topological, komlos1996topological} and has since found several important applications, see e.g., \cite{kim2017proof,alon2023essentially, liu2023solution,bucic2022towards} and the nice survey by Letzter~\cite{letzter2024sublinear}). It is easy to prove that sublinear expanders can be found in essentially any graph. Moreover, recent work by Buci\'c and Montgomery~\cite{bucic2022towards} (using some ideas from Tomon~\cite{MR4695962}) shows that random vertex subsets in sublinear expanders with polylogarithmic average degree are likely to inherit some expansion properties. Thus, if the graph $H$ above was also a sublinear expander (see Subsection~\ref{sec:sublin} for the notion of expansion we use), we could change our partition of $V(H)$ to $R,V_1,\ldots,V_t$, before proceeding as outlined above to obtain edge-disjoint linear forests $F^i$, $i \in [k]$, with vertex set $\cup_{j=1}^t V_j$ and then using vertices from $R$ to connect each linear forest $F^i$ into a cycle $C^i$.

In order to successfully combine these approaches by regularity and by expansion, we will need to answer the following meta-question.

\begin{qn*}
Given any graph with polylogarithmic average degree, can we find a dense subgraph $H$ which is extremely close to being regular and which has some expansion properties?
\end{qn*}

Despite regularity and sublinear expansion being key properties that have been used to prove various results about sparse graphs, there has been no good answer to this meta-question so far in the literature. Perhaps the closest is by Dragani\'c, Methuku, Munh\'a Correia and Sudakov~\cite{draganic2023cycles}, who recently used sublinear expanders which are almost-regular, but the 100-almost-regularity they achieve is far too weak for the proposed outline above.
A key contribution of this paper is to address this issue by developing a novel near-regularisation technique for a graph $G$ that is (say) already 100-almost-regular, allowing us to find a subgraph $H$ of $G$ which is not only extremely close to being regular but also contains a large random set of vertices $A\subset V(G)$, and which contains a good proportion of the density of $G$. In conjunction with our careful development of the techniques in \cite{bucic2022towards} (for transferring some expansion properties into the random subset $A$), this would allow us to carry out the outline above to find $k$ edge-disjoint cycles with almost the same vertex set.

This new regularisation technique has the potential to be an important tool in this area, and variants of this tool allow it to be effective even in graphs with only (large) constant average degree. For example, in our upcoming paper \cite{CJMMregular}, building on this technique, we develop powerful methods for finding regular (or nearly-regular) subgraphs of large degree. Variants of this lemma will also be used in \cite{MMPSpathdecomp}. For more on this, see the concluding remarks in \Cref{sec:concluding}. We now state the variant we prove for this paper.

\begin{lemma}[Regularisation lemma]
\label{thm:effnearreg2}
For every $\lambda> 1$, there is some $C\geq 1$ such that the following holds. Let $G$ be an $n$-vertex graph in which every vertex has degree between $d$ and $\lambda d$. Let $A\subset V(G)$ be chosen by including each vertex of $G$ independently at random with probability $1/C$.
Then, with probability $1-o(1)$, $G$ contains a $(d'\pm 10^5 \lambda^5\log n)$-nearly-regular subgraph~$H$ with $d/C\leq d'\leq d$ and $A\subset V(H)$. 
\end{lemma}

Due to its independent interest, we give a sketch of the proof of \Cref{thm:effnearreg2} in
Subsection~\ref{sec:newreg} before giving a complete proof in \Cref{sec:regularizing lemma}.
To reiterate, what makes this lemma particularly effective for our application (and potentially others) is that it ensures that a large random subset $A$ of vertices of the original graph $G$ is contained in the nearly-regular subgraph $H$ (with high probability), allowing us to exploit the expansion properties of $G$. To illustrate this, suppose that $G$ is an expander. Then it is quite possible that the nearly-regular subgraph $H$ we find in it is not an expander (making $H$ potentially useless for making the path connections required for building our edge-disjoint cycles with the same vertex set).
However, fortunately, since $A$ is a random subset of the expander $G$, the lemma still allows us to use the expansion properties of $G$ to connect vertices of $H$  using paths in $G$ through the random subset~$A$.

\subsection{The full problem: incorporating absorption}\label{sec:fullproblem}
In order to find $k$ edge-disjoint cycles with exactly the same vertex set, we will use absorption (the powerful method first codified by R\"odl, Ruci\'nski and Szemer\'edi~\cite{rodl2006dirac}). 
A natural absorption strategy applied to the outline in Subsection~\ref{sec:halfproblem} would be to find a path $Q^i$ (for each $i \in [k]$) which can absorb any set $L^i\subset R$ of vertices that we do not use for connecting the paths in the linear forest $F^i$. Here, this would mean that there is a path $Q^i_*$ with the same endvertices as $Q^i$ which has vertex set $V(Q^i)\cup L^i$. To simplify matters, imagine that we could even find $Q^i$ as one of the paths within $F^i$. Then, substituting $Q^i$ by $Q^i_*$ in the cycle $C^i$ in the outline (given in Subsection~\ref{sec:halfproblem}) will give a cycle with vertex set exactly $V(F^i)\cup R=V(H)$ for all $i \in [k]$.

It is relatively easy to find a suitable absorbing path $Q^i$ using expansion properties; however, we cannot guarantee that it is one of the paths within $F^i$, so the difficulty here is to find $Q^i$ in a way that it combines well with our approach requiring near-regularity. More precisely, suppose that the final vertex set of the cycles $C^i$ that we aim to find is $V\cup R$ for some $V\subset V(H)\setminus R$. Then, note that, for the strategy outlined in Subsection~\ref{sec:halfproblem} to succeed, we would like to have that $H[V\setminus V(Q^i)]$ is very nearly regular, for each $i\in [k]$. (Indeed, if this near-regularity does hold, then we can find a spanning linear forest $F^i$ with few components in $H[V\setminus V(Q^i)]$, and we can connect the endpoints of the paths in $F^i$ and the endpoints of $Q^i$ using paths through $R$ to form a cycle, which can absorb the unused vertices of $R$, yielding a cycle $C^i$ with vertex set $V\cup R$ for all $i \in [k]$.) 
We will achieve this near-regularity by constructing the absorbing paths $Q^i$ by concatenating $s$ randomly selected short absorbing paths, each able to absorb just one pair of vertices. More precisely, for some $s$, and vertex set $X=\{x_1,\ldots,x_{s+1}\}$, each path $Q^i$ will consist of the concatenation of $s$ short absorbing paths, each joining $x_j$ and $x_{j+1}$ for some $j \in [s]$, and with $(\log n)^{O(1)}$ vertices. To randomise these $s$ short absorbing paths, we do the following. Using some set $U$ disjoint from $X$, we will construct for each $j\in [s]$, $k$ short absorbing paths joining  $x_{j}$ and $x_{j+1}$ (using vertices from $U$)
with the absorption property we want (see \Cref{def:absorberov} and the discussion after it) so that, across all $j \in [s]$ the short absorbing paths found are vertex disjoint outside of $X$. For each $j\in [s]$, independently at random, we will assign the $k$ short absorbing $x_{j}$-$x_{j+1}$ paths to the $k$ cycles we are constructing, and then, for each $i\in [k]$, concatenate these $s$ short absorbing paths to form the absorbing path $Q^i$ assigned to the $i$-th cycle.

Let $U_\mathrm{abs} \subseteq U$ be the set of all the vertices outside of $X$ used in the absorbing paths $Q^i$ over all $i \in [k]$. Then, because of the randomised construction of $Q^i$, for each vertex $v$, the expected degree (in $H$) of $v$ into $U_\textrm{abs}\setminus V(Q^i)$ is $(1-\frac{1}{k})d_H(v,U_{\mathrm{abs}})$. Take $W\subset U\setminus U_{\mathrm{abs}}$ by including each vertex uniformly at random with probability $(1-\frac{1}{k})$. Then, for every vertex $v$, the expected degree of $v$ into $W\cup (U_\mathrm{abs}\setminus V(Q^i))$ is equal to $(1-\frac{1}{k})d_H(v,U)$. As the short absorbing paths are, indeed, short, there will be sufficient concentration that with high probability every vertex $v$ will have around $(1-\frac{1}{k})d_H(v,U)$ neighbours in $W\cup (U_\mathrm{abs}\setminus V(Q^i))$.
Thus, if all the vertex degrees in $H$ are approximately the same into $U$ (which will hold because $H$ is very nearly regular and $U$ is a random subset), then with high probability, we will have that $H[W\cup (U_\mathrm{abs}\setminus V(Q^i))]$ is approximately regular for each $i\in [k]$. Hence (by the discussion in the previous paragraph, letting $V = W\cup U_\mathrm{abs}\cup X$) we can then carry out the strategy outlined in Subsection~\ref{sec:halfproblem} to obtain $k$ edge-disjoint cycles with vertex set $W \cup U_\mathrm{abs}\cup X \cup R$, as desired.

Until now, we have ignored one complication which we now address. Our graph may be bipartite, and, indeed, as is standard we may as well assume it is so.
Therefore, an absorbing path will only be able to `absorb' a balanced set of vertices, that is, with the same number of vertices in each side of the bipartition. To be more precise, we use the following notion of absorbers for our `short absorbing paths'.

\begin{defn}[Absorber]
\label{def:absorberov}
Given distinct vertices $a$, $b$, $y$, $z$, and a set $S$ of vertices with $a,b,y,z\not \in S$, an absorber for the pair $a, b$ is the union of two paths joining $y$ and $z$, one with internal vertex set $S$, and another with internal vertex set $S \cup \{a, b\}$. We say that $S$ is the interior of the absorber, and $y$ and $z$ are the endpoints of the absorber.
\end{defn}

In the above outline, for every $i \in [k]$, we would like to partition the `leftover set' $L^i\subset R$ into balanced vertex pairs, before absorbing the vertex pairs into the absorbing path $Q^i$. It will not be hard to make sure that $L^i$ is balanced, but we cannot take an absorber for every possible balanced pair of vertices, as there are too many possible pairs. Therefore, we use a `robustly matchable' bipartite graph (introduced by Montgomery~\cite{montgomery2019spanning}) as
an auxiliary graph to tell us for which pairs of vertices we need to construct absorbers. In more detail, we take random balanced sets $R_1, R_2$, and let $R = R_1 \cup R_2$. Let $K$ be a robustly matchable bipartite graph with $s$ edges, with small maximum degree and with vertex set $R_1\cup R_2$. We will ensure that each of the $s$ short absorbing paths used to construct the absorbing path $Q^i$ is capable of absorbing a unique pair $\{a,b\}$ with $ab\in E(K)$. Then, once we have used (a balanced set of) vertices from only $R_1$ to connect the paths in the linear forest $F^i$, $K$ will have a perfect matching even after removing the vertices in $R_1$ which have been used (since it is `robustly matchable'). We can then absorb all of the vertex pairs in this matching into $Q^i$ so that every vertex in $R = R_1\cup R_2$ is contained in the cycle $C^i$.

\subsection{Detailed outline}\label{sec:detailedoutline}
We will now make precise the ideas introduced in the previous two subsections and give a detailed overview of our proof. Our aim is to find $k$ edge-disjoint cycles on the same vertex set in a graph $G$ with sufficiently many edges.

\begin{enumerate}[label=(\alph*)]
\itemsep 0.5 mm
    \item \label{overview:expander} Find a bipartite subgraph $G'$ of $G$ which is a suitable expander such that $G'$ is an $n$-vertex graph with average degree at least $(\log n)^{C}$ (where $C$ is a sufficiently large constant), where the vertex degrees differ by at most a factor of 18 (using Lemma~\ref{cor:findexpander}).

    \item \label{overview:regsub} Using Lemma~\ref{thm:effnearreg2}, find a subgraph $H\subset G'$ that is very nearly regular, almost as dense as $G'$, and which contains a random vertex subset $A$.

    \item \label{overview:nicesets} Take pairwise disjoint random balanced subsets $R_1, R_2, X, U$ in $V(H)$ of appropriate sizes, where, as $A\cap R_1\subset R_1$ and $A\cap U\subset U$, with high probability we will be able to connect vertex pairs in $H$ using paths (in the expander $G'$) through $R_1$ and $U$.
    Let $X \coloneqq \{x_1, \ldots, x_{s+1}\}$. We think of $R \coloneqq R_1 \cup R_2$ as a random `reservoir'. 
\end{enumerate}

\noindent \textbf{Constructing absorbers}
\begin{enumerate}[label=(\alph*), resume]
\itemsep 0.5 mm
    \item \label{overview:RMBG} Take a collection $K \coloneqq \{p_1, \ldots, p_s\}$ of $\Theta(|R|)$ many balanced pairs in $R$ such that $K$ is `robustly matchable' in the following sense: for every balanced set $R'_1\subseteq R_1$ with $|R'_1|\le |R_1|/2$, there is a perfect matching in $K[(R_1\cup R_2)\setminus R'_1]$. (Here, with a slight abuse of notation, we denote by $K$ the graph whose vertex set is $R$ and whose edge set is~$K$.)

    \item \label{overview:absorbers} Using the expansion properties of $G'$, for every $p_j \in K$, construct $k$ absorbers $\Abs_j^i$, $i \in [k]$, with endpoints $x_j$ and $x_{j+1}$ whose interiors $S_j^i$ are pairwise disjoint subsets of $U$ of size at most $(\log n)^{12}$.
    Here, each of the absorbers $\Abs_j^i$ contains two paths joining $x_j$ and $x_{j+1}$, one with internal vertex set $S_j^i$, and another with internal vertex set $S_j^i \cup p_j$ (i.e., $\Abs_j^i$ contains a path joining $x_j$ and $x_{j+1}$ which can `absorb' the pair $p_j$). Let $U_{\textrm{abs}}$ be the set of vertices of $U$ used in these absorbers. 

    \item \label{absorbing path}  For every $p_j \in K$, independently and randomly assign each of the $k$ absorbers $\Abs_j^i$, $i \in [k]$, to a different cycle that we aim to construct. For every $i \in [k]$, let $U_{\textrm{abs}}^i$ be the (balanced) subset of $U_{\textrm{abs}}$ consisting of those vertices which are used in the absorbers assigned to the $i$-th cycle.

    Observe that the union of absorbers assigned to the $i$-th cycle contains a path $Q^i$  (joining $x_1$ and $x_{s+1}$) with vertex set $U_{\textrm{abs}}^i \cup X$ such that any pair in $K$ can be `absorbed' into $Q^i$.

\item \label{balancing set} Choose a balanced set $W\subset U\setminus U_{\textrm{abs}}$ of size $(1-\frac{1}{k})|U\setminus U_{\textrm{abs}}|$ uniformly at random.

\end{enumerate}

\noindent \textbf{Constructing {\boldmath $k$} edge-disjoint linear forests}

\begin{enumerate}[label=(\alph*), resume]
\itemsep 0.5 mm
    \item \label{definingVi} For every $i\in [k]$, let $V^i= W\cup(U_{\textrm{abs}}\setminus U_{\textrm{abs}}^i)$. (Note that $V^i$ is approximately distributed like a random subset of $U$ where each element of $U$ is chosen with probability $1 - \frac{1}{k}$.) Let $V^i_1,\dots,V^i_t$ be a random partition of $V^i$ into $t$ balanced sets of equal size (where $t$ is chosen appropriately).

    \item \label{find matchings} For every $i\in [k]$, find nearly perfect matchings $M^i_j$ in $H[V^i_j,V^i_{j+1}]$ for all $j\in [t-1]$ such that all these matchings are pairwise edge-disjoint, and edge-disjoint from the absorbers. (This can be done since $H[V^i_j,V^i_{j+1}]$ is very nearly regular due to the fact that $H[U]$ is very nearly regular and $V^i$ is distributed like a random subset of $U$.) For every $i\in [k]$, let $F^i$ be the linear forest obtained by taking the union of the matchings $M^i_j$ for $j\in [t-1]$.
\end{enumerate}

\pagebreak

\noindent \textbf{Extending the linear forests to {\boldmath $k$} edge-disjoint cycles on the same vertex set}

\begin{figure}[h]
    \centering
\includegraphics[width=0.8 \textwidth]{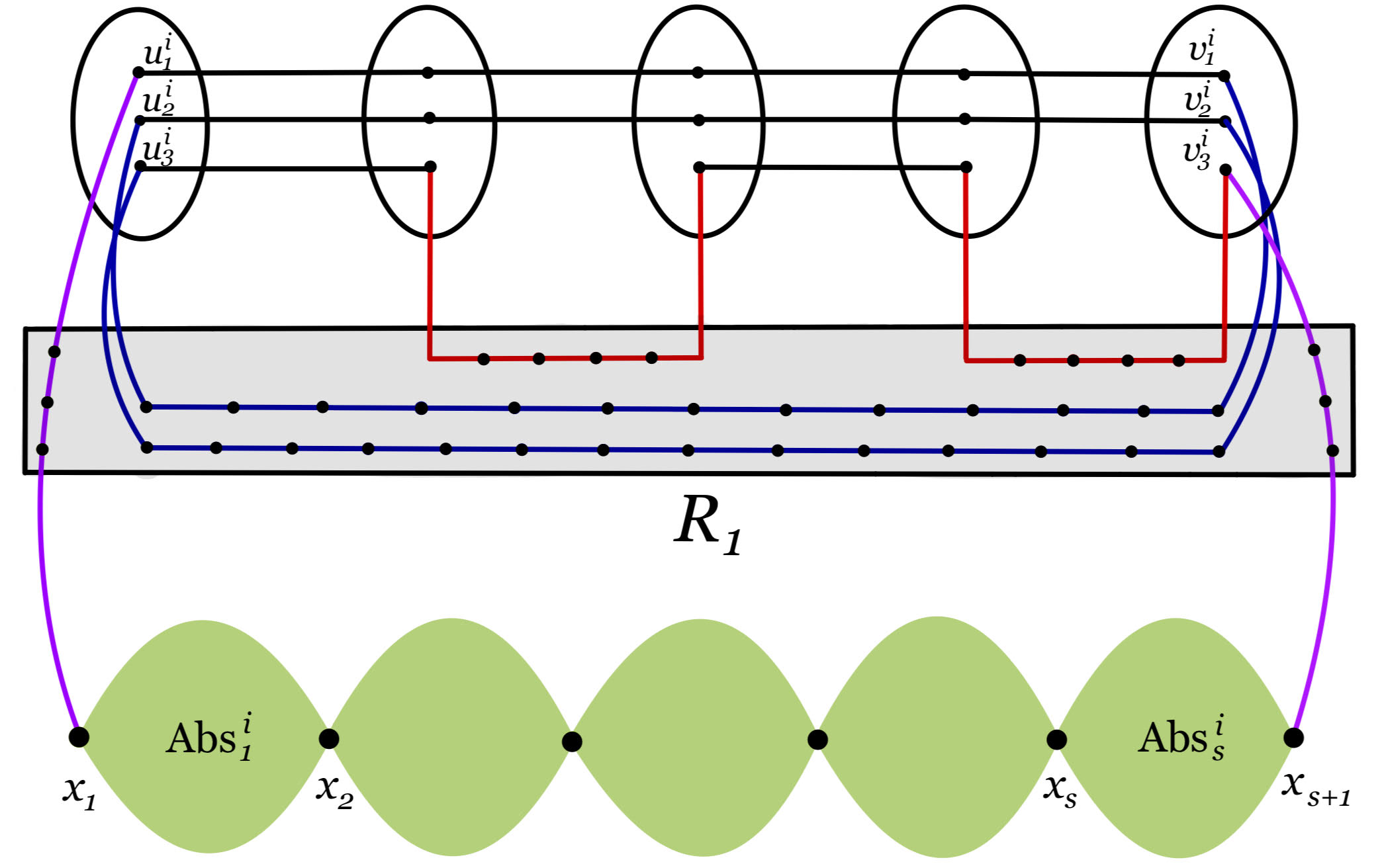}
    \caption{Constructing the $i$-th cycle $C^i$.}
\label{fig:constructingcycle}
\end{figure}

\begin{enumerate}[label=(\alph*), resume]
\itemsep 0.5 mm

\item \label{ov:connectingpaths} Using the expansion properties of $G'$, connect the paths in $F^i$ and the endpoints $x_1$ and $x_{s+1}$ of the `absorbing' path $Q^i$  through $R_1$ (see Figure~\ref{fig:constructingcycle}, where the connecting paths are the red, blue and purple paths). Ensure, in addition, that all these connecting paths are internally vertex-disjoint and edge-disjoint from the absorbers.

\item \label{path containing matchings} These connecting paths together with
$F^i$ yield a path $P^i$ joining $x_1$ and $x_{s+1}$ with internal vertex set $V^i \cup R^i_1 = W\cup (U_{\textrm{abs}}\setminus U_{\textrm{abs}}^i) \cup  R^i_1$, where $R^i_1 \subseteq R_1$ is the set of vertices in $R_1$ used in these connecting paths. We will ensure that $|R^i_1| \le |R_1|/2$. Note that the paths $P^i$ for $i \in [k]$ are pairwise edge-disjoint since all of the connecting paths are internally vertex-disjoint, and the linear forests $F^i$ (given by \ref{find matchings}) are pairwise edge-disjoint. Furthermore, each $P^i$ is edge-disjoint from the absorbers too (as both the linear forests $F^i$ and the connecting paths through $R_1$ are).

\item \label{ov:robust}Since $K$ is `robustly matchable',  \ref{overview:RMBG} ensures that there is a perfect matching $K^i$ in $K[(R_1\setminus R^i_1) \cup R_2]$ for every $i \in [k]$. By `absorbing' all of the pairs in $K^i$ into the path $Q^i$ (given by \ref{absorbing path}), we obtain a path $Q^i_*$ with endpoints $x_1$ and $x_{s+1}$ and with vertex set $V(Q^i_*) = U^i_{\textrm{abs}} \cup X \cup (R_1\setminus R^i_1) \cup R_2$. The path $Q^i_*$ together with the path $P^i$ (given by \ref{path containing matchings}) yields a cycle $C^i$ with vertex set $W\cup U_{\textrm{abs}} \cup X \cup R$ for every $i \in [k]$. Hence, we found $k$ edge-disjoint cycles on the same vertex set, as desired.
\end{enumerate}

\subsection{Finding very nearly regular subgraphs}\label{sec:newreg}

Let us now give a sketch of the proof of \Cref{thm:effnearreg2}.  We will obtain the desired nearly-regular subgraph in $G$ by finding a sequence $G=:G_0\supset G_1\supset G_2\supset \dots$ of subgraphs such that each $G_{i+1}$ is a `bit more regular' than $G_i$ but has similar average degree. We terminate the process once $G_i$ is sufficiently regular and let $H$ be this final subgraph. Let us explain how, given a graph $G_i$ whose degrees are between $d_i$ and $(1+\gamma)d_i$, we find a subgraph $G_{i+1}\subset G_i$ with better regularity properties. Let $U_{L}=\{v\in V(G_i): d_{G_i}(v)\leq (1+\gamma/2)d_i\}$ be the set of `low-degree vertices' in $G_i$, and let $U_{H}=\{v\in V(G_i): d_{G_i}(v)> (1+\gamma/2)d_i\}$ be the set of `high-degree vertices' in $G_i$. Let $\eps=\gamma/100$ and let $G_{i+1}$ be the subgraph of $G_i$ obtained by
\begin{itemize}[noitemsep]
    \item  deleting edges within $U_H$ independently at random with probability $2\eps-\eps^2$,
    \item deleting edges from $U_{L}$ to $U_{H}$ independently at random with probability $\eps$, and
    \item deleting vertices in $U_{L}$  independently at random with probability $\eps$.
\end{itemize}
It is easy to see that if $v\in U_L$, then $$\mathbb{E}[d_{G_{i+1}}(v) \hspace{1mm} | \hspace{1mm} v\in V(G_{i+1})]=(1-\eps)d_{G_i}(v),$$
whereas if $v\in U_H$, then
$$\mathbb{E}[d_{G_{i+1}}(v)]=(1-\eps)^2 d_{G_i}(v).$$
This means that, on average, the degrees of high-degree vertices drop faster than the degrees of low-degree vertices. Hence, (provided that we have suitable concentration), $G_{i+1}$ will have slightly better regularity properties than $G_i$. This can be continued as long as degrees differ by more than an additive term of roughly $\log n$, as we then have enough concentration to take a union bound over all vertices in $G_i$. Thus, we eventually arrive at a graph that is extremely close to being regular. The other crucial property, namely that $V(H)$ contains a large random subset $A$ of $V(G)$, will follow from the observation that in each step we keep each vertex $v\in V(G_i)$ with probability at least $1-\eps$, so even after taking many steps of this process, vertices in $G$ `survive' with at least a positive constant probability.

\section{Key lemmas and tools}
\label{sec:key lemmas}

One of the key results we will use in the proof of \Cref{thm:main} is our regularisation lemma (\Cref{thm:effnearreg2}), discussed in the proof sketch (Section~\ref{sec:overview}). In the following subsections we present the other main tools that we will need in our proof.

\subsection{Sublinear expansion}
\label{sec:sublin}

We use the following definition of an expander graph from \cite{bucic2022towards} (see \cite{bucic2022towards} for more details of the development of this style of expansion, including related definitions developed independently by Haslegrave, Kim, and Liu~\cite{haslegrave2022extremal} and by Sudakov and Tomon~\cite{sudakov2022extremal}).
\begin{restatable}{defn}{defnsublinear}\label{defn:robust-sublinear-expansion}
An $n$-vertex graph $G$ is an $(\eps,s)$-expander if, for every $U\subseteq V(G)$ and $F\subseteq E(G)$ with $1\le |U|\leq \frac23n$ and $|F|\leq s|U|$, we have
\begin{equation}
|N_{G-F}(U)|\geq \frac{\eps|U|}{(\log n)^2}.\label{eqn:expands}
\end{equation}
\end{restatable}

The following lemma finds an almost-regular (bipartite) expander in any graph while losing at most a logarithmic factor in the average degree. In our proof of Theorem \ref{thm:main}, we use this lemma to pass to such an expander.

\begin{restatable}{lemma}{findexpander}\label{cor:findexpander}
    Let $n$ be a sufficiently large integer and let $0< \eps < 2^{-3}$. Let $G$ be an $n$-vertex graph with $d(G)\geq (\log n)^4$. Then, $G$ contains an $18$-almost-regular bipartite subgraph $G'$ with $d(G')\geq \frac{d(G)}{400\log n}$ which is an $(\eps,s)$-expander for some $s\geq \frac{d(G')}{(\log |V(G')|)^2}$. 
\end{restatable}

The proof of \Cref{cor:findexpander} is by now standard. Following, for example, \cite{draganic2023cycles} we first find within $G$ a $6$-almost-regular bipartite subgraph with average degree at least $d(G)/100\log n$ (using a result proved in \cite{ARS+17} and \cite{bucic2020nearly}), before finding within this an expander with large minimum degree using the variation in \cite{bucic2022towards} on the standard proofs to find sublinear expanders. For completeness, we include a full proof of Lemma~\ref{cor:findexpander} in the appendix.

\subsection{Connecting pairs of vertices using vertex-disjoint paths through a random set}

In the proof of \Cref{thm:main}, we will often connect pairs of vertices with internally vertex-disjoint paths through a given set. The following definition describes a set through which we can simultaneously connect many pairs of vertices. The condition on the pairs we connect is a little intuitive as it restricts the number of neighbours (with multiplicity) some vertices can have among the vertices in the pairs. This condition will be convenient for us to use as we will be working in almost-regular subgraphs $G$, where we will have $D\approx d(G)^{1-\eps}$ for some small $\eps>0$.

\begin{defn} [Connecting set] \label{def:nice}
    Let $G$ be an $n$-vertex graph. We say that a set $V\subset V(G)$ is $D$-connecting (in $G$) if for every sequence $x_1,\dots,x_r,y_1,\dots,y_r$ of (not necessarily distinct) vertices outside of $V$ such that each vertex in $V$ has at most $D$ neighbours in the multiset $\{x_1,\dots,x_r,y_1,\dots,y_r\}$, there is a collection of internally vertex-disjoint paths $P_{x_i y_i}$, $i \in [r]$, such that, for each $i \in [r]$, $P_{x_i y_i}$ is an $x_i-y_i$ path (in $G$) through $V$ of length at most $(\log n)^6$.
\end{defn}

The following lemma shows that random subsets of vertices of (a sufficiently regular) expander are connecting sets.

\begin{lemma} \label{lem:randomsetisnice}
    Let $G$ be an $n$-vertex $(\eps,s)$-expander with minimum degree $\delta$ and maximum degree $\Delta$, where $2^{-90}< \eps< 1$ and $s\ge p^{-4}(\log n)^{30}$ for some $\frac{100\log n}{\delta}\leq p\leq 1$. Let $V$ be a $p$-random subset of $V(G)$. Then, with probability $1-o(1)$, $V$ is $D$-connecting for $D=\frac{p^9 s \delta}{\Delta (\log n)^{73}}$.
\end{lemma}
We prove this lemma in \Cref{sec: sublinear expander}. Its proof is an adaptation of arguments from \cite{bucic2022towards} together with some new ideas. At the beginning of \Cref{sec: sublinear expander}, we provide a detailed overview of how our result and argument differ from those in \cite{bucic2022towards}.

\subsection{Absorbers}
\label{subsec:absorbers}
For convenience, we recall here our definition for an absorber.

\begin{defn}[Absorber]
\label{def:absorber}
Given distinct vertices $a$, $b$, $y$, $z$, and a set $S$ of vertices with $a,b,y,z\not \in S$, an absorber for the pair $a, b$ is the union of two paths joining $y$ and $z$, one with internal vertex set $S$, and another with internal vertex set $S \cup \{a, b\}$, neither of which contains any edges between the vertices $a$, $b$, $y$ and $z$. We say that $S$ is the interior of the absorber, and $y$ and $z$ are the endpoints of the absorber.
\end{defn}

\begin{figure}[h]
    \centering
\includegraphics[width=0.75 \textwidth]{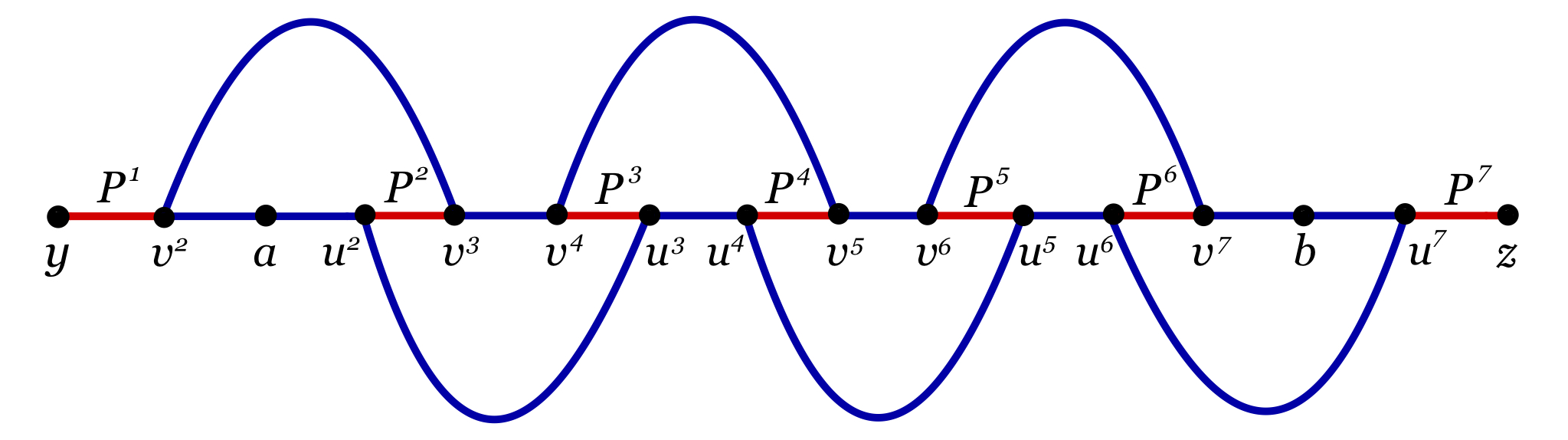}
    \caption{The absorber for the pair $a, b$ when $t = 8$. The edges are coloured blue, and the paths $P^1, P^2, \ldots, P^7$ are coloured red.}
    \label{fig:absorber}
\end{figure}

Note that, compared to \Cref{def:absorberov}, we have added the technical condition that there is no edge between $a$, $b$, $y$  and $z$ in the absorbers.
We will now construct an absorber for a given pair of vertices $a, b$ (see Figure~\ref{fig:absorber}), using a natural construction that appears in, for example, \cite{montgomery2019spanning}.

\textbf{An absorber for the pair $a,b$.} Let $t$ be an even integer. Let $u^1, u^2, \ldots ,u^t$ and $v^1, v^2, \ldots, v^t$ be two paths such that $u^1 = v^1 = a$, $u^t = v^t = b$ and $u^i \ne v^j$ for all $1 < i,j < t$. Take two new vertices, $y$ and $z$. Let $P^1$ be a path from $y$ to $v^2$, let $P^{t-1}$ be a path from $u^{t-1}$ to $z$, and for $2\leq i\leq t-2$ let $P^i$ be a path joining $u^i$ and $v^{i+1}$. Suppose, moreover, that the paths $P^1, P^2, \ldots, P^{t-1}$ are pairwise vertex-disjoint and do not contain $a$ or $b$. Let $S = (\cup_{i=1}^{t-1} V(P^i))\setminus \{a,b,y,z\}$.

Consider the following two paths. 
\begin{itemize}
\item $y, P^1,v^2,
a, u^2,P^2, v^3, v^4, P^3, \ldots,u^{t-3},u^{t-2},P^{t-2},v^{t-1},b, u^{t-1}, P^{t-1}, z$.
\item $y, P^1, v^2, v^3, P^2, u^2, u^3, P^3, v^4, v^5, P^4, u^4,\ldots, u^{t-3}, P^{t-3}, v^{t-2}, v^{t-1}, P^{t-2}, u^{t-2}, u^{t-1}, P^{t-1}, z$.
\end{itemize}

Note that the first path has internal vertex set $S \cup \{a, b\}$, and the second path has internal vertex set $S$ (see \Cref{fig:absorber}). Hence, the union of these paths is an absorber for $a, b$ with endpoints $y$ and $z$.

The next lemma allows us to construct absorbers for given pairs of vertices by first connecting the pairs using paths through a connecting set $U_1$, and then making further connections using paths through another connecting set $U_2$.
(The existence of such connecting sets will come from \Cref{lem:randomsetisnice}.)

\begin{lemma} \label{lem:finddisjointabsorbers}
Let $G$ be a bipartite graph with parts $A$ and $B$, and let $U_1, U_2\subseteq V(G)$ be disjoint sets such that $U_i$ is $D_i$-connecting for each $i\in [2]$. Let $y_1,\ldots,y_r,z_1,\ldots,z_r\in V(G)$, $a_1,\ldots,a_r\in A$ and $b_1,\ldots,b_r\in B$ be (not necessarily distinct) vertices outside of $U_1\cup U_2$ such that $y_i,z_i,a_i$ and $b_i$ are distinct for each $i\in [r]$,
\begin{enumerate}[\rm(i)] 
    \item\label{findabsorbersa} every vertex in $U_1$ has at most $\frac{D_1}{(\log n)^6}$ neighbours in the multiset $\{a_1,\ldots,a_r,b_1,\ldots,b_r\}$, and
    \item \label{findabsorbersb} every vertex in $U_2$ has at most $D_2$ neighbours in the multiset $U_1\cup \{y_1,\ldots,y_r,z_1,\ldots,z_r\}$.
\end{enumerate}
Then, for each $i\in [r]$ there is an absorber $\Abs_i$ for the pair $a_i,b_i$ with endpoints $y_i$ and $z_i$ such that the interiors of these absorbers are pairwise disjoint subsets of $U_1\cup U_2$ of size at most $(\log n)^{12}$. 
\end{lemma}

\begin{proof} 
Consider the sequence $\sigma$ obtained by replacing every element of the sequence $a_1,\dots,a_r,b_1,\dots,b_r$ with $(\log n)^6$ copies of that element. Then,  by \ref{findabsorbersa}  of the lemma, every vertex in $U_1$ has at most $D_1$ neighbours in the multiset consisting of the elements of $\sigma$. Therefore, since $U_1$ is $D_1$-connecting (see \Cref{def:nice}), we can find $(\log n)^6$ many $a_i$-$b_i$ paths of length at most $(\log n)^6$ with internal vertices in $U_1$ for every $i\in [r]$ such that all these paths are pairwise internally vertex-disjoint. 
Thus, for every $i\in [r]$, we can find two paths $Q_i, Q'_i$ between $a_i$ and $b_i$ such that the lengths of $Q_i$ and $Q'_i$ are the same, the paths $Q_i, Q'_i$ have their internal vertices in $U_1$, and all these paths are pairwise internally vertex-disjoint. For every $i \in [r]$, let $Q_i = u_i^1 u_i^2 \ldots u_i^{t_i}$ and $Q'_i = v_i^1 v_i^2 \ldots v_i^{t_i}$, where $u_i^1 = v_i^1 = a_i$ and $u_i^{t_i} = v_i^{t_i} = b_i$ (see e.g., the two blue paths between $a$ and $b$ in \Cref{fig:absorber}). 
Note that $t_i$ must be even because $a_i$ and $b_i$ are on opposite sides of the bipartition $A\cup B$. 
Next, since $U_2$ is $D_2$-connecting, using \ref{findabsorbersb} of the lemma, we can find paths $P_i^1, P_i^2, \ldots, P_i^{t_i-1}$ of length at most $(\log n)^6$ joining the pairs $(y_i,v_i^2), (u_i^2,v_i^3), (u_i^3,v_i^4), \ldots, (u_i^{t_i-2},v_i^{t_i-1}), (u_i^{t_i-1},z_i)$, respectively, for every $i\in [r]$ 
(see e.g., the red paths in \Cref{fig:absorber}) such that 
all these paths are pairwise internally vertex-disjoint, and 
the internal vertices of all these paths are in $U_2$. For each $i\in [r]$, let $S_i = \cup_{j=1}^{t_i-1}V(P_i^j) \setminus \{a_i,b_i,y_i,z_i\}$. Then, for each $i\in [r]$, by the discussion before Figure~\ref{fig:absorber}, there exists an absorber $\Abs_i$ for the pair $a_i,b_i$ with interior $S_i$ and endpoints $y_i$ and $z_i$. As the sets $S_i$, $i\in [r]$, are pairwise disjoint subsets of $U_1\cup U_2$ of size at most $(\log n)^{6} \cdot (\log n)^{6} = (\log n)^{12}$ by construction, this completes the proof of the lemma.
\end{proof}

We use the following lemma to construct our absorbers in the proof of Theorem~\ref{thm:main} (see Subsection~\ref{proofsec:absorbers}).

\begin{lemma} \label{finding absorber chain}
Let $k$ be a positive integer and let $G$ be a bipartite graph with parts $A$ and $B$. Let $R,X,U_1, U_2\subseteq V(G)$ be disjoint sets such that $U_i$ is $D_i$-connecting for $i\in [2]$. Let $K\subseteq (R\cap A)\times (R\cap B)$ be a set of vertex pairs such that no vertex of $R$ appears in more than 200 pairs in $K$. Assume that $|X|=|K|+1$ and write $X=\{x_1,\ldots,x_{s+1}\}$ and $K= \{p_1,\ldots,p_s\}$. Suppose that
\begin{enumerate}[\rm(i)]
    \item \label{findabsorberchaina} every vertex $v\in U_1$ satisfies $d_G(v,R)\le \frac{D_1}{200k(\log n)^6}$, and
    \item \label{findabsorberchainb} every vertex $v\in U_2$ satisfies $d_G(v,U_1\cup X)\le \frac{D_2}{2k}$.
\end{enumerate}
Then, for all $i\in [k]$ and $j\in [s]$ there is an absorber $\Abs_j^i$ for the pair $p_j$ with endpoints $x_j$ and $x_{j+1}$ such that the interiors of these absorbers (for all $i\in [k]$ and $j\in [s]$) are pairwise disjoint subsets of $U_1\cup U_2$ of size at most $(\log n)^{12}$. 
\end{lemma}

\begin{proof}
Let $r=ks$. For $j\in [s]$, let $p_j = (p_j^1,p_j^2)$. 
For every $i\in [r]$, find the unique $j\in [s]$ such that $s$ divides $i-j$, and set $a_i = p_j^1$, $b_i = p_j^2$, $y_i=x_j$ and $z_i=x_{j+1}$. Since each vertex in $R$ appears in at most $200$ pairs $(p_j^1,p_j^2)\in K$, it follows that each vertex in $R$ appears at most $200k$ times in the multiset $\{a_1,\dots,a_r,b_1,\dots,b_r\}$. Thus, by \ref{findabsorberchaina} of the lemma, every vertex in $U_1$ has at most $\frac{D_1}{(\log n)^6}$ neighbours in the multiset $\{a_1,\dots,a_r,b_1,\dots,b_r\}$. Moreover, each vertex in $X$ appears at most $2k$ times in the multiset $\{y_1,\dots,y_r,z_1,\dots,z_r\}$ (and every element of the multiset is in $X$). Hence, by \ref{findabsorberchainb} of the lemma, every vertex in $U_2$ has at most $D_2$ neighbours in the multiset $U_1\cup \{y_1,\dots,y_r,z_1,\dots,z_r\}$.
Thus, we can apply \Cref{lem:finddisjointabsorbers} to find the desired absorbers.
\end{proof}

\subsection{Robustly matchable bipartite graphs}

We use robustly matchable bipartite graphs as auxiliary graphs to tell us which pairs of vertices to construct absorbers for (see Subsection~\ref{proofsec:absorbers} for more details). These graphs were introduced in \cite{montgomery2019spanning} to find spanning trees in random graphs.

\begin{lemma} [Lemma~10.7 in \cite{montgomery2019spanning}] \label{robustly-matchable-bipartite-graphs}
For every sufficiently large $m$, there exists a bipartite graph $H$ with maximum degree at most $100$ with vertex classes $X$ and $Y\cup Z$, with $|X|=3m$ and $|Y|=|Z|=2m$ such that the following holds. For every $Z'\subseteq Z$ with $|Z'|=m$, there is a perfect matching in $H$ between $X$ and $Y\cup Z'$.
\end{lemma}

We deduce the following bipartite analogue of the above lemma, which is essential for our proof.

\begin{lemma} \label{lem:bipartite graph absorbing absorber}
For every sufficiently large $m$, there exists a bipartite graph $H$ with an odd number of edges, with maximum degree at most $102$, and with vertex classes $A_1\cup A_2$ and $B_1\cup B_2$, with $|A_1|=|B_1|=2m$ and $|A_2|=|B_2|=5m$ such that the following holds. For every $A'_1\subseteq A_1$ and $B'_1\subseteq B_1$ with $|A'_1|=|B'_1|\ge m$, there is a perfect matching in $H$ between $A_1'\cup A_2$ and $B_1'\cup B_2$.
\end{lemma}

\begin{proof}
Consider a bipartite graph $H$ satisfying the property in \Cref{robustly-matchable-bipartite-graphs} with $m$, and take two vertex-disjoint copies of it, $H_1$ and $H_2$, where $H_i$ has vertex classes $X_i$ and $Y_i\cup Z_i$. Define $A_1 = Z_1$, $B_1 = Z_2$, $A_2 = X_2 \cup Y_1$ and $B_2 = X_1 \cup Y_2$. 
Consider a perfect matching $M$ between $A_1$ and $B_1$. Consider the graph $H$ that is the union of $H_1$, $H_2$, and $M$. 
If this graph $H$ contains an even number of edges, then we add an arbitrary edge to it so the number of edges is odd. 
We next show that this graph satisfies the other desired properties. Clearly, $H$ has maximum degree at most $102$ and $|A_1|=|B_1|=2m$ and $|A_2|=|B_2|=5m$. Now, suppose that $A'_1\subseteq A_1$ and $B'_1\subseteq B_1$ with $|A'_1|=|B'_1|\ge m$. Then, first, find a matching $M' \subseteq M$ between $A'_1$ and $B'_1$ of size $|A'_1| - m$. Let $A''_1 = A'_1 \setminus V(M')$ and $B''_1 = B'_1 \setminus V(M')$. Since $|A''_1| = |B''_1| = m$, by the properties of $H_i$, we can find a perfect matching $M_1$ between $X_1$ and $Y_1 \cup A''_1$ and another perfect matching $M_2$ between $X_2$ and $Y_2 \cup B''_1$. By taking the union of $M'$, $M_1$, and $M_2$, we get the desired perfect matching between $A_1'\cup A_2$ and $B_1'\cup B_2$.
\end{proof}

\subsection{Concentration inequalities}\label{sec:conc}
We will use the following basic version of the Chernoff bound for the binomial random variable (see, for example, \cite{alon2016probabilistic}).

\begin{lemma}[Chernoff bound] \label{chernoff}
Let $n$ be an integer and $0\le \delta,p \le 1$. If $X \sim \bin(n,p),$ then, setting $\mu=\mathbb{E} [X] = np,$ we have
$$\P(X>(1+\delta) \mu) \le e^{-\delta^2\mu/2}\quad\quad\quad \text{ and }\quad\quad\quad \P(X<(1-\delta) \mu) \le e^{-\delta^2\mu/3}.$$
\end{lemma}

We will also often use the following well-known martingale concentration result (see Chapter 7 of \cite{alon2016probabilistic}).

\begin{lemma}\label{lem:mcd}
Suppose that $X:\prod_{i=1}^N\Omega_i\to \R$ is $k$-Lipschitz. Then, for each $t>0$,
\[
\mathbb{P}(|X-\mathbb{E} [X]|>t)\leq 2\exp\left(\frac{-t^2}{2k^2N}\right).
\]
\end{lemma}

\section{Finding a nearly-regular subgraph containing a random vertex subset}
\label{sec:regularizing lemma}

In this section we prove our regularisation lemma (\Cref{thm:effnearreg2}) following the sketch in Subsection~\ref{sec:newreg}. At the heart of its proof is the following lemma, which shows, roughly speaking, that in any graph $G$ we can find a subgraph with similar average degree to $G$ but with slightly better regularity properties, and which in addition contains a very large random subset of the vertices.

\begin{lemma}\label{lem:onestephighprob} Let $\eps,\gamma>0$ with $\eps\leq 1/100$ and $\gamma\geq 10\eps$. 
Let $n$ and $d$ be such that $\eps d\geq 10^4\log n$.
Let $G$ be a graph on at most $n$ vertices such that $d\leq d_G(v)\leq (1+\gamma)d$ for each $v\in V(G)$. Let $A\subset V(G)$ be chosen by including each vertex of $G$ independently at random with probability $1-\eps$. Then, with probability at least $1-\frac{1}{n^2}$, $G$ contains a subgraph $G'$ such that $A\subset V(G')$ and, for some $(1-2\eps)d\leq d'\leq d$, we have $d'\leq d_{G'}(v)\leq (1-\frac{\eps}{2})(1+\gamma)d'$ for each $v\in V(G')$.
\end{lemma}

\begin{proof} Let $U_L=\{v\in V(G):d_G(v)\leq (1+\gamma/2)d\}$ and $U_H=\{v\in V(G):d_G(v)> (1+\gamma/2)d\}$ be the set of low and high degree vertices in $G$, respectively.

Let $G'$ be the random subgraph of $G$ obtained by
\begin{enumerate}[\rm(i)] 
    \item  deleting edges within $U_H$ independently at random with probability $2\eps-\eps^2$, \label{item:between dense}
    \item deleting edges from $U_{L}$ to $U_{H}$ independently at random with probability $\eps$, and \label{item:sparse and dense}
    \item deleting vertices in $U_{L}$  independently at random with probability $\eps$. \label{item:delete vertices}
\end{enumerate}

Note that, equivalently, $V(G')$ can be generated by taking $V(G')=(A\cap U_L)\cup U_H$. Hence, we may assume that $A\subset V(G')$.

Now note that if $v\in U_L$ and $uv\in E(G)$, then conditional on $v\in V(G')$, the probability that $uv\in E(G')$ is precisely $1-\eps$. Indeed, if $u\in U_L$, then (conditional on $v\in V(G')$) $uv\in E(G')$ holds if and only if $u$ did not get deleted by \ref{item:delete vertices}, which has probability $1-\eps$, and if $u\in U_H$, then (conditional on $v\in V(G')$) $uv\in E(G')$ holds if and only if the edge $uv$ did not get deleted by \ref{item:sparse and dense}, which has probability $1-\eps$.

Moreover, if $v\in U_H$ and $uv\in E(G)$, then the probability that $uv\in E(G')$ is precisely $(1-\eps)^2=1-2\eps+\eps^2$. Indeed, if $u\in U_L$, then $uv\in E(G')$ holds if and only if $u$ did not get deleted by \ref{item:delete vertices} and the edge $uv$ did not get deleted by \ref{item:sparse and dense}, which has probability $(1-\eps)^2$, whereas if $u\in U_H$, then $uv\in E(G')$ holds if and only if the edge $uv$ did not get deleted by \ref{item:between dense}, which has probability $1-(2\eps-\eps^2)$.

Hence, if $v\in U_L$, then conditional on $v\in V(G')$ the distribution of $d_{G'}(v)$ is $\bin(d_G(v),1-\eps)$, and if $v\in U_H$, then the distribution of $d_{G'}(v)$ is $\bin(d_G(v),(1-\eps)^2)$.

It follows from the Chernoff bound (Lemma~\ref{chernoff}) that if $v\in U_L$, then, as $d_G(v)\geq d$, 
$$\mathbb{P}\left(d_{G'}(v)< \Big(1-\frac{5}{4}\eps\Big)d \hspace{1mm} | \hspace{1mm} v\in V(G')\right)\leq \mathbb{P}\left(\bin(d,1-\eps)< \Big(1-\frac{5}{4}\eps\Big)d\right)=\mathbb{P}\left(\bin(d,\eps)> \frac{5}{4}\eps d\right)\leq e^{-(1/4)^2 \eps d/2}\leq \frac{1}{n^4}.$$

Moreover, if $v\in U_H$, then
\begin{align*}
    \mathbb{P}\left(d_{G'}(v)> \Big(1-\frac{7}{4}\eps\Big)(1+\gamma)d \right)
    &\leq \mathbb{P}\left(\bin((1+\gamma)d,(1-\eps)^2)> \Big(1-\frac{7}{4}\eps\Big)(1+\gamma)d\right) \\
    &= \mathbb{P}\left(\bin((1+\gamma)d,2\eps-\eps^2)< \frac{7}{4}\eps (1+\gamma) d\right) \\
    &\leq \mathbb{P}\left(\bin\Big((1+\gamma)d,\frac{15}{8}\eps\Big)< \frac{7}{4}\eps (1+\gamma) d\right)\leq e^{-(1/15)^2 \eps d/3}\leq \frac{1}{n^4},
\end{align*}
where the penultimate inequality follows from the Chernoff bound (Lemma~\ref{chernoff} applied with $\delta=1/15$), using that $\frac{7}{4}\eps(1+\gamma)d=(1-\frac{1}{15})\cdot \mathbb{E}[\bin((1+\gamma)d,\frac{15}{8}\eps)]$ and $\eps d\leq \mathbb{E}[\bin((1+\gamma)d,\frac{15}{8}\eps)]$. Furthermore, if $v\in U_H$, then
\begin{align*}
    \mathbb{P}\left(d_{G'}(v)< \Big(1-\frac{5}{4}\eps\Big)d\right)
    &\leq \mathbb{P}\left(d_{G'}(v)< d\right) \leq \mathbb{P}\left(\bin((1+\gamma/2)d,(1-\eps)^2)< d\right)\leq \mathbb{P}\left(\bin((1+\gamma/2)d,1-2\eps)< d\right) \\
    &= \mathbb{P}\left(\bin((1+\gamma/2)d,2\eps)> \gamma d/2\right) \leq e^{-\eps d}\leq \frac{1}{n^4},
\end{align*}
where the penultimate inequality follows from the Chernoff bound since $2\eps d\leq \E[\bin((1+\gamma/2)d,2\eps)]\leq \gamma d/4$.

Note also that, if $v\in U_L\cap V(G')$, then  $d_{G'}(v)\leq (1+ \frac{\gamma}{2})d\leq (1-\frac{7}{4}\eps)(1+\gamma)d$ (with probability $1$). Thus, for every $v\in V(G)$, conditional on $v\in V(G')$, the probability that $(1-\frac{5}{4}\eps)d\leq d_{G'}(v)\leq (1-\frac{7}{4}\eps)(1+\gamma)d$ is at least $1-\frac{1}{n^3}$.

Hence, by the union bound, the probability that $(1-\frac{5}{4}\eps)d\leq d_{G'}(v)\leq (1-\frac{7}{4}\eps)(1+\gamma)d$ for each $v\in V(G')$ is at least $1-\frac{1}{n^2}$. Setting $d'=(1-\frac{5}{4}\eps)d$, we have
\[
d'\leq d_{G'}(v)\leq \left(1-\frac{7}{4}\eps\right)(1+\gamma)d=\frac{1-\frac{7}{4}\eps}{1-\frac{5}{4}\eps}(1+\gamma)d'\leq \left(1-\frac{\eps}{2}\right)(1+\gamma)d'
\]
for each $v\in V(G')$, as required.
\end{proof}

\begin{proof}[ of \Cref{thm:effnearreg2}]
    Let $C$ be sufficiently large (in terms of $\lambda$) to support our argument. If $d< C\log n$, then we may take $d'=d/C$ and let $H$ be the empty graph on vertex set $V(G)$, so let us assume that $d\geq C\log n$. Let $\eps=\frac{10^4\lambda^5\log n}{d}$. Since $C$ is sufficiently large and $d\geq C\log n$, we have $\eps\leq 1/100$. Let $k$ be the smallest non-negative integer such that $\lambda(1-\eps/2)^k\leq 1+10\eps$. Since $\lambda(1-\eps/2)^n\leq \lambda e^{-\eps n/2}\leq \lambda e^{-\eps d/2}<1$, we have $k\leq n$. Moreover, clearly $\lambda(1-\eps/2)^k\geq 1$, so $(1-\eps/2)^k\geq 1/\lambda$.

    We will prove by induction on $i$ that
    
    $(\ast)$ for each $0\leq i\leq k$, if $A_i$ is a $(1-\eps)^i$-random subset of $V(G)$, then with probability at least $1-\frac{i}{n^2}$, $G$ contains a subgraph $G_i$ such that $A_i\subset V(G_i)$ and for some $(1-2\eps)^i d\leq d_i\leq d$, we have $d_i\leq d_{G_i}(v)\leq \lambda (1-\eps/2)^i d_i$ for each $v\in V(G_i)$.

    Assuming that $(\ast)$ holds, we have in particular (taking $i=k$) that if $A_k$ is a $(1-\eps)^k$-random subset of $V(G)$, then with probability at least $1-\frac{k}{n^2}=1-o(1)$, $G$ contains a subgraph $G_k$ such that $A_k\subset V(G_k)$ and for some $(1-2\eps)^k d\leq d_k\leq d$, we have $d_k\leq d_{G_k}(v)\leq \lambda (1-\eps/2)^k d_k$ for each $v\in V(G_k)$. Then we have $d_k\leq d_{G_k}(v)\leq  (1+10\eps)d_k \le d_k +10\eps d = d_k+10^5 \lambda^5 \log n$ for each $v\in V(G_k)$, so $G_k$ is $(d_k\pm 10^5\lambda^5 \log n)$-nearly-regular. Moreover, since $(1-\eps/2)^{k}\geq 1/\lambda$, it follows that $(1-2\eps)^k\geq 1/\lambda^5$, hence $d_k\geq d/\lambda^5\geq d/C$. Finally, since $(1-\eps)^k\geq 1/\lambda^5\geq 1/C$, $A_k$ (and hence $V(G_k)$) contains a suitable random subset of $V(G)$, so $H=G_k$ and $d'=d_k$ satisfy the conditions of the lemma.

    It remains to prove $(\ast)$. Note that $(\ast)$ trivially holds for $i=0$. Assume now that it holds for some $0\leq i<k$. Let $A_i\subset V(G)$ have the property that $G$ contains a subgraph $G_i$ such that $A_i\subset V(G_i)$ and, for some $(1-2\eps)^i d\leq d_i\leq d$, we have $d_i\leq d_{G_i}(v)\leq \lambda (1-\eps/2)^i d_i$ for each $v\in V(G_i)$. By the definition of $k$, we have $\lambda (1-\eps/2)^i\geq 1+10\eps$. Moreover, since $(1 - \eps/2)^i \ge 1/\lambda$,
    we have $d_i\geq (1-2 \eps)^i d \ge d/\lambda^5$ as before, so $\eps d_i\geq 10^4\log n$. Now let $R_{i+1}$ be a $(1-\eps)$-random subset of $V(G_i)$, and let $A_{i+1}=R_{i+1}\cap A_i$. Then, by applying \Cref{lem:onestephighprob} (with $G_i$ and $R_{i+1}$ playing the roles of $G$ and $A$ respectively), 
    with probability at least $1-\frac{1}{n^2}$,  $G_i$ contains a subgraph $G_{i+1}$ such that $A_{i+1} \subset R_{i+1}\subset V(G_{i+1})$ and, for some $(1-2\eps)d_i\leq d_{i+1}\leq d_i$, we have $d_{i+1}\leq d_{G_{i+1}}(v)\leq \lambda (1-\eps/2)^{i+1} d_{i+1}$ for all $v\in V(G_{i+1})$. Moreover, $A_{i+1}$ is a $(1-\eps)$-random subset of $A_i$. Hence, if $A_i$ is a $(1-\eps)^i$-random subset of $V(G)$, then $A_{i+1}$ is a $(1-\eps)^{i+1}$-random subset of $V(G)$. This completes the induction step and the proof of the lemma.
\end{proof}

\section{Further auxiliary lemmas}
\label{sec:auxiliary lemmas}

We now give the results that we use to find the linear forests (in Subsection~\ref{sec:matchings}) and well-behaved vertex subsets (in Subsections~\ref{sec:degreesto} and~\ref{sec:findingvertex}) in the proof of Theorem~\ref{thm:main}.

\subsection{Matchings covering most vertices}\label{sec:matchings}

We will use the following lemma to construct our linear forests in the proof of Theorem~\ref{thm:main} (see Subsection~\ref{subsec:linearforestsproof}) by taking the union of almost-perfect matchings that are chosen randomly. Crucially, the lemma allows us to appropriately bound the degrees of vertices into the set $Y$ of `unmatched' vertices (which is essential for connecting these vertices using paths through $R_1$ in \ref{ov:connectingpaths} in Subsection~\ref{sec:detailedoutline} of the proof sketch).

\begin{lemma} \label{lem:matchings with little leftover}
    Let $G$ be a graph and let $V_1,\dots,V_t$ be disjoint subsets of $V(G)$ such that, for each $j\in [t-1]$, the degree of every vertex in $G[V_j,V_{j+1}]$ is between $\delta$ and $\Delta$. Let $T_1,\dots,T_n$ be subsets of $V(G)$ such that, for each $i\in [n]$ and $j\in [t]$, we have $|T_i\cap V_j|\leq r$. Assume that $1-\delta/\Delta\leq t^{-1/2}\log n$. Then, there exist matchings $M_j$ for each $j\in [t-1]$ in $G[V_j,V_{j+1}]$ such that, if $Y$ consists of all those vertices $y\in V_1\cup V_t$ which belong to neither $M_1$ nor $M_{t-1}$ and all those vertices $y\in V_j$ with $2\leq j\leq t-1$ which do not belong to both $M_{j-1}$ and $M_j$, then $|Y|\le 10|\bigcup_{i=1}^t V_i|t^{-1/2}\log n$ and $|Y\cap T_i|\leq 10rt^{1/2}\log n$ for all $i\in [n]$.
\end{lemma}

\begin{proof}
    It is well known that every bipartite graph with maximum degree at most $\Delta$ has a proper edge-colouring with $\Delta$ colours. Let us take such a colouring of $G[V_j,V_{j+1}]$ for each $j\in [t-1]$. The colour classes give a collection of at most $\Delta$ matchings in $G[V_j,V_{j+1}]$ which partition the edge set of $G[V_j,V_{j+1}]$. For each $j\in [t-1]$ independently, let $M_j$ be the matching defined by a uniformly random colour class in $G[V_j,V_{j+1}]$. For each vertex $u\in G[V_j,V_{j+1}]$, the probability that $u$ belongs to $M_j$ is at least $d(u)/\Delta$ where $d(u)$ is the degree of $u$ in the graph $G[V_j,V_{j+1}]$. Hence, the probability that $u$ does not belong to $M_j$ is at most $1-\delta/\Delta$. It follows that, for each $i\in [t]$, we have $\mathbb{E}[|Y\cap T_i|]\leq tr \cdot 2(1-\delta/\Delta)\leq 2rt^{1/2}\log n$. Moreover, the choice for $M_j$ changes the value of $|Y\cap T_i|$ by at most $2r$ (since $|T_i\cap V_j|,|T_i\cap V_{j+1}|\leq r$). Hence, by \Cref{lem:mcd}, $$\mathbb{P}\left(|Y\cap T_i|\geq 10rt^{1/2}\log n\right)\leq \mathbb{P}\left(|Y\cap T_i|-\mathbb{E}[|Y\cap T_i|]\geq 5rt^{1/2}\log n\right)\leq 2\exp\left(\frac{-(5rt^{1/2}\log n)^2}{2(2r)^2 t}\right)\leq \frac{1}{2n}.$$
    Thus, by the union bound over all $i\in [n]$, with probability at least $1/2$ we have $|Y\cap T_i|\leq 10rt^{1/2}\log n$ for all $i\in [n]$.

    Furthermore, since the probability that any given vertex $u\in V(G)$ is in $Y$ is at most $2(1-\delta/\Delta)\leq 2t^{-1/2}\log n$, it follows that $\mathbb{E}[|Y|]\leq 2|\bigcup_{i=1}^t V_i|t^{-1/2}\log n$. Hence, by Markov's inequality, the probability that $|Y|\geq 10|\bigcup_{i=1}^t V_i|t^{-1/2}\log n$ is at most $1/5$. Thus, with positive probability, the matchings $M_j$ for $j \in [t-1]$ satisfy the desired properties.
\end{proof}

\subsection{Degrees to a random subset of vertices}\label{sec:degreesto}
For \ref{find matchings} in Subsection~\ref{sec:detailedoutline} (of the proof sketch), we want to preserve some near-regularity conditions while taking a random vertex subset, for which we use the following lemma.
\begin{lemma} \label{lem:simple random subgraph}
    Let $G$ be an $n$-vertex graph. Let $V\subseteq V(G)$. Let $d\ge (\log n)^{15}$, let $0\le d' \le d$ and let $1\le t\le d$ be an integer. 
    Let $U$ be a uniformly random subset of $V$ of size $\lfloor \frac{1}{t}|V| \rfloor$. 
    \begin{enumerate}
        \item[$\mathrm{(a)}$] If $v\in V(G)$ satisfies $d_G(v,V)\le d+ d'$, then $d_G(v,U)\le \frac{1}{t}d + d' + d^{2/3}$ with probability $1-n^{-\omega(1)}$. 
        \item[$\mathrm{(b)}$] If $v\in V(G)$ satisfies $d_G(v,V)\ge d- d'$, then $d_G(v,U)\ge \frac{1}{t}d - d' - d^{2/3}$ with probability $1-n^{-\omega(1)}$. 
    \end{enumerate}
\end{lemma}

We will deduce this from the following more general technical lemma which allows, in addition to the uniformly random subset $U$, some further random subsets in which vertices can appear with some dependencies (as we need for picking our absorbers randomly in \ref{absorbing path} in Subsection~\ref{sec:detailedoutline} (of the proof sketch)). 
\Cref{lem:simple random subgraph} easily follows from the following result by setting $V = V'$, $\tau = 1$, and $s=0$.

\begin{lemma} \label{lem:random subgraph}
    Let $G$ be an $n$-vertex graph. Let $V'\subseteq V(G)$. Let $d\ge (\log n)^{15}$, $0\le d' \le d$ and let $1\leq t\leq d$ and $s\geq 0$ be integers. Let $\tau=1$ or $\tau=t-1$. 
    Let  $V^i_j$, $i\in [t]$ and $j\in [s]$, be pairwise disjoint sets in $V(G)\setminus V'$ satisfying $|V^i_j|\le  d^{1/8}$. Let $V = V'\cup \left(\bigcup_{i\in [t], j\in [s]} V^i_j\right)$.  
    For every $j\in [s]$, let $S_j$ be a uniformly random subset of $[t]$ of size $\tau$. 
    Let $U$ be a uniformly random subset of $V'$ of size $\lfloor \frac{\tau}{t}|V'| \rfloor$. 
    Let $Z$ denote the random set $U\cup \left(\bigcup_{j\in [s], i\in S_j} V^i_j\right)$.
    \begin{enumerate}
        \item[$\mathrm{(1)}$] If $v\in V(G)$ satisfies $d_G(v,V)\le d+ d'$, then $d_G(v,Z)\le \frac{\tau}{t}d + d' + d^{2/3}$ with probability $1-n^{-\omega(1)}$. 
        \item[$\mathrm{(2)}$] If $v\in V(G)$ satisfies $d_G(v,V)\ge d- d'$, then $d_G(v,Z)\ge \frac{\tau}{t}d - d' - d^{2/3}$ with probability $1-n^{-\omega(1)}$. 
    \end{enumerate}
\end{lemma}

\begin{proof}
    To avoid unimportant technicalities, assume that $|V'|$ is divisible by $t$. 
    Note that it suffices to prove the lemma in the special case where $V'=\emptyset$. Indeed, in the general case, let us take an enumeration $\big\{u^i_j : i \in [t], s+1 \le j \le s+\frac{|V'|}{t} \big\}$ of the elements of the set $V'$ uniformly at random and define new sets $V^i_j \coloneqq \{ u^i_j\}$ for $i \in [t]$, $s+1 \le j \le s+\frac{|V'|}{t}$. Then we can apply the lemma in the $V'=\emptyset$ case (with $s$ replaced by $s+\frac{|V'|}{t}$ and with the sets $V_j^i$ for all $i\in [t]$ and $j\in \big[s+\frac{|V'|}{t}\big]$), since the random subset $Z$ has the same distribution as before.
    Thus, for the rest of the proof we assume that $V'=\emptyset$.

    We have $\mathbb{E}[d_G(v,Z)]=\frac{\tau}{t}d_G(v,V)$. Since $\tau=1$ or $\tau=t-1$ and $|V_j^i|\leq d^{1/8}$ for each $i\in [t]$ and $j\in [s]$, changing the choice of $S_j$ (for some $j\in [s]$) can change the value of $d_G(v,Z)$ by at most $d^{1/8}$. 
    Note that the random variable $d_G(v,Z)$ depends on at most $d_G(v,V) \le d + d' \le 2d$ choices of $S_j$, $j\in [s]$, namely only those $j\in [s]$ for which $N_G(v)\cap \left(\bigcup_{i\in [t]} V^i_j\right) \neq \emptyset$. Thus, for each $v\in V(G)$, we can apply \Cref{lem:mcd} with some $N\leq 2d$ for both (1) and (2), to show that the corresponding event does not occur with probability at most 
    $$\mathbb{P}\left(\big| d_G(v,Z)-\mathbb{E}[d_G(v,Z)] \big| \geq d^{2/3}\right)
    \leq 2\exp\left(\frac{-(d^{2/3})^2}{2(d^{1/8})^2 \cdot 2d}\right)= n^{-\omega(1)},$$
    completing the proof.
\end{proof}

\subsection{Finding vertex sets with good degree control and good connectivity properties}\label{sec:findingvertex}

The following lemma will be used to find a suitable subgraph $G'$ and connecting sets $R_1,R_2,X,U_1,U_2\subset V(G')$ in the proof of Theorem~\ref{thm:main} (see Subsection~\ref{proofsubsec:nicesets}). The lemma ensures that the degrees of vertices to the sets $R_1,R_2,X,U_1$ in $G$ are bounded appropriately. More crucially, it also ensures that $G'[U_1 \cup U_2] = G'[U]$ is extremely close to being regular (which is essential for constructing nearly perfect matchings, as explained in \ref{find matchings} in Subsection~\ref{sec:detailedoutline} of the proof sketch). 

\begin{lemma} \label{lem:findsets}
    There is an absolute constant $C$ such that the following holds. Let $t$ be a fixed positive integer, let $n$ be sufficiently large and let $G$ be a bipartite $n$-vertex $18$-almost-regular $(\eps,s)$-expander with average degree $d=d(G)\geq (\log n)^C$, where $s\geq d/(\log n)^2$ and $2^{-90}< \eps< 1$. Let $m_1,\dots,m_t$ be positive integers greater than $n/d^{1/10}$ and smaller than $n/(t\log n)$. Then, there exists a graph $G'\subset G$ and pairwise disjoint subsets $S_1,\dots,S_t\subset V(G')$ such that, for each $i\in [t]$, we have that 
    
    \begin{itemize}
        \item[$\mathrm{(i)}$] $S_i$ is balanced (with respect to the bipartition of $G$),
        \item[$\mathrm{(ii)}$] $|S_i|=2m_i$,
        \item[$\mathrm{(iii)}$] $S_i$ is $D_i$-connecting in $G$ for $D_i=\frac{p_i^9 d}{(\log n)^{76} }$, where $p_i=\frac{m_i}{n}$, 
        
        \item[$\mathrm{(iv)}$] there exists some $d_i$ such that $\frac{dm_i}{Cn} \le d_i \le \frac{Cdm_i}{n}$, and $d_i-d_i^{2/3}\leq d_{G'}(v,S_i)\leq d_i+d_i^{2/3}$ holds for every $v\in V(G')$ and
        \item[$\mathrm{(v)}$] $d_G(v,S_i)\leq \frac{Cdm_i}{n}$ holds for every $v\in V(G)$.
    \end{itemize}
\end{lemma}

\begin{proof}
We will take $C$ to be large enough so that the following arguments hold when $d\geq (\log n)^C$.
By \Cref{thm:effnearreg2} (with $\lambda = 18$), there exists some constant $C_0$ such that if $A$ is a $(1/C_0)$-random subset of $V(G)$, then with probability $1 - o(1)$, $G$ contains a $(d'\pm C_0\log n)$-nearly-regular subgraph $G'$, where $d/C_0\leq d'\leq d$ and $A\subset V(G')$. Let $n'=|V(G')|$. Note that, since $A \subset V(G')$, with probability $1-o(1)$ we have $n'\geq n/(2C_0)$. Furthermore, for each $i\in [t]$, if $B_i\subset A$ is chosen by including each vertex independently at random with probability $m_i/n$, then $B_i$ is a $(m_i/nC_0)$-random subset of $V(G)$. Thus, by \Cref{lem:randomsetisnice} and as $G$ is $18$-almost-regular, $m_i > n/d^{1/10}$, and $s \ge \frac{d}{(\log n)^2}$, we have that $B_i$ is $D_i$-connecting in $G$ with probability $1-o(1)$.

Therefore, in total, we can conclude that we can choose some $G'\subset G$ with the following properties. 
\begin{enumerate}[label = (\alph*)]
\item $G'$ is $(d'\pm C_0\log n)$-nearly-regular (for some $d/C_0\leq d'\leq d$).
\item $n'\geq n/(2C_0)$ (where $n'=|V(G')|)$.
\item For each $i\in [t]$ and any $p\geq m_i/n$, a $p$-random subset of $V(G')$ is $D_i$-connecting in $G$ with probability $1-o(1)$. 
\end{enumerate}

We will now choose disjoint sets $R_i,T_i\subset V(G')$, $i\in [t]$, randomly and show that, with high probability, any sets $S_i$ with $R_i\subset S_i\subset R_i\cup T_i$ for each $i\in [t]$, will satisfy each of (iii)--(v), and, furthermore, that we can pick such sets $S_i$ such that (i) and (ii) hold.

For this, randomly take disjoint sets $R_i,T_i$, $i\in [t]$, in $V(G')$, by, for each vertex $v\in V(G')$, independently placing $v$ in each $R_i$ with probability $\frac{2m_i}{n'}-\frac{4C_0m_i\log n}{d'n'}-\frac{m_i^{3/5}}{n'}$ and in each $T_i$ with probability $\frac{10C_0m_i\log n}{d'n'}+\frac{2m_i^{3/5}}{n'}$. Note that this is possible as, for each $i\in [t]$, $\frac{4C_0m_i\log n}{d'n'}+\frac{m_i^{3/5}}{n'}\leq \frac{2m_i}{n'}$, and
\[
\sum_{i\in [t]}\bigg(\frac{2m_i}{n'}-\frac{4C_0m_i\log n}{d'n'}-\frac{m_i^{3/5}}{n'}+\frac{10C_0m_i\log n}{d'n'}+\frac{2m_i^{3/5}}{n'}\bigg)\leq \frac{2n}{n'\log n}+\frac{10C_0n}{d'n'}+\frac{2n}{n'\log n}\leq 1,
\]
as $m_i\leq n/(t\log n)$ for all $i\in [t]$, $n'\geq n/(2C_0)$ and $d'\geq d/C_0\geq (\log n)^C/C_0$.
Furthermore, for each $i\in [t]$, choose $S_i$ arbitrarily such that $R_i\subset S_i\subset R_i\cup T_i$, and, if possible, satisfying (i) and (ii). 
We will now show that each of (i)--(v) hold with high probability, so that, as $n$ is large, we can take the required sets $S_i$, $i\in [t]$. Note first that, by (c) we have that (iii) holds with probability $1-o(1)$ since $R_i$ is a $p$-random subset of $V(G')$ for $p=\frac{2m_i}{n'}-\frac{4C_0 m_i \log n}{d' n'}-\frac{m_i^{3/5}}{n'}\geq \frac{m_i}{n}$.
  
 Let $X\cup Y$ be the bipartition of $G'$ inherited from the unique bipartition of $G$, labelling so that $|X|\geq |Y|$. Then, as $G'$ is $(d'\pm C_0\log n)$-nearly-regular, we have 
 \[
 |Y|\geq \frac{1-\frac{C_0\log n}{d'}}{1+\frac{C_0\log n}{d'}}\cdot |X|\geq \bigg(1-2\cdot \frac{C_0\log n}{d'}\bigg)\frac{n'}{2},\;\;\;\text{ and }\;\;\; |X|= n'- |Y|\leq \bigg(1+2\cdot \frac{C_0\log n}{d'}\bigg)\frac{n'}{2},
 \]
 so that  $\frac{n'}{2}- \frac{C_0n'\log n}{d'}\leq |X|,|Y|\leq \frac{n'}{2}+ \frac{C_0n'\log n}{d'}$. Since $|X|\leq \frac{n'}{2}+ \frac{C_0n'\log n}{d'}$, for each $i\in [t]$ we have
    $$\mathbb{E}[|R_i\cap X|]\leq \frac{n'}{2}\cdot \left(\frac{2m_i}{n'}-\frac{4C_0m_i\log n}{d'n'}-\frac{m_i^{3/5}}{n'}\right)+\frac{C_0n'\log n}{d'}\cdot \frac{2m_i}{n'}= m_i-\frac{m_i^{3/5}}{2}.$$
As $m_i\geq n/d^{1/10}\geq n^{9/10}$, it follows by the Chernoff bound that with probability at least $1 - \frac{1}{n}$, we have $|R_i\cap X|\leq m_i$. On the other hand, $|X|\geq \frac{n'}{2}- \frac{C_0n'\log n}{d'}$ and the probability that some $v\in V(G')$ belongs to $R_i\cup T_i$ is $\frac{2m_i}{n'}+\frac{6C_0m_i\log n}{d'n'}+ \frac{m_i^{3/5}}{n'}$, so 
    
    \begin{align*}
        \mathbb{E}[|(R_i\cup T_i)\cap X|]
        &\geq \left(\frac{n'}{2}-\frac{C_0n'\log n}{d'}\right)\cdot \left(\frac{2m_i}{n'}+\frac{6C_0m_i\log n}{d'n'}+\frac{m_i^{3/5}}{n'}\right) \\
        &\geq \frac{n'}{2}\cdot \left(\frac{2m_i}{n'}+\frac{6C_0m_i\log n}{d'n'}+\frac{m_i^{3/5}}{n'}\right)-\frac{C_0n'\log n}{d'}\cdot \frac{3m_i}{n'} = m_i+\frac{m_i^{3/5}}{2}.
    \end{align*}
    
Hence, by the Chernoff bound, we have, with probability at least $1 - \frac{1}{n}$, that $|(R_i\cup T_i)\cap X|\geq m_i$. By symmetry, we also have $|R_i\cap Y|\leq m_i\leq |(R_i\cup T_i)\cap Y|$ with probability at least $1 -\frac{1}{n}$. Therefore, by a union bound, it follows that for all $i \in [t]$, we have $|R_i\cap X|\leq m_i\leq |(R_i\cup T_i)\cap X|$ and $|R_i\cap Y|\leq m_i\leq |(R_i\cup T_i)\cap Y|$ with probability at least $1 - \frac{2t}{n}=1-o(1)$, so with high probability we will have selected the $S_i$, $i\in [t]$, so that (i) and (ii) hold.

As $G$ is $18$-almost-regular with average degree $d$, each vertex $v\in V(G)$ has $d_G(v)\leq 18d$. 
Note that, then, as the probability that a given vertex $u\in V(G')$ belongs to $R_i\cup T_i$ is at most $\frac{3m_i}{n'}$, for any $i\in [t]$ and $v\in V(G)$, we have $\mathbb{E}[d_G(v,R_i\cup T_i)]\leq d_G(v)\cdot \frac{3m_i}{n'}\leq 18d\cdot \frac{3m_i}{n'}\leq \frac{108C_0dm_i}{n}\leq \frac{Cdm_i}{2n}$ (provided that $C$ is sufficiently large). Hence, by the Chernoff bound, the probability that $d_G(v,S_i)>\frac{Cdm_i}{n}$ is at most $\frac{1}{n^2}$. It follows by the union bound that with probability $1-o(1)$, $\mathrm{(v)}$ of the lemma is satisfied.

For each $i\in [t]$, let $d_i=2d'm_i/n'$. Note that since $d/C_0\leq d'\leq d$, $n/(2C_0)\leq n' \leq n$, and $C$ can be chosen larger than $4 C_0$, we have $\frac{d m_i}{C n} \le \frac{2d m_i}{C_0 n} \le d_i \le \frac{4C_0 d m_i}{n} \le \frac{C d m_i}{n}$. Now let $v\in V(G')$. Since $d_{G'}(v)\geq d'-C_0\log n$, we have
    $$\mathbb{E}[d_{G'}(v,R_i)]\geq d'\cdot \left(\frac{2m_i}{n'}-\frac{4C_0m_i\log n}{d'n'}-\frac{m_i^{3/5}}{n'} \right)-C_0\log n\cdot \frac{2m_i}{n'}= d_i-6C_0\log n\cdot \frac{m_i}{n'}-d'\frac{m_i^{3/5}}{n'}.$$
Since $6C_0\log n\cdot m_i/n'+d'm_i^{3/5}/n'\leq (\log n)^2 +d_i/(2m_i^{2/5})\leq d_i^{2/3}/2$, we have $\mathbb{E}[d_{G'}(v,R_i)]\geq d_i-d_i^{2/3}/2$, therefore, by the Chernoff bound, with probability at least $1-\frac{1}{n^2}$, we have $d_{G'}(v,R_i)\geq d_i-d_i^{2/3}$. Furthermore, since $d_{G'}(v)\leq d'+C_0\log n$ and $d'\geq d/C_0\geq (\log n)^C/C_0$, we have
    \begin{align*}
        \mathbb{E}[d_{G'}(v,R_i\cup T_i)]
        &\leq (d'+C_0\log n)\cdot \left(\frac{2m_i}{n'}+\frac{6C_0m_i\log n}{d'n'}+\frac{m_i^{3/5}}{n'}\right)\leq d'(1+d_i^{-1/2})\cdot \left(\frac{2m_i}{n'}\right)(1+10d_i^{-2/5}) \\
        &=d_i(1+d_i^{-1/2})(1+10d_i^{-2/5})\leq d_i+d_i^{2/3}/2.
    \end{align*}
Hence, with probability at least $1-\frac{1}{n^2}$, we have $d_{G'}(v,R_i\cup T_i)\leq d_i+d_i^{2/3}$. Thus, by the union bound, with probability at least $1 - \frac{2t}{n}$, we have $d_{G'}(v,R_i)\geq d_i-d_i^{2/3}$ and $d_{G'}(v,R_i\cup T_i)\leq d_i+d_i^{2/3}$ for all $v\in V(G')$ and $i \in [t]$. In this case, $d_i-d_i^{2/3}\leq d_{G'}(v,S_i)\leq d_i+d_i^{2/3}$ for all $v\in V(G')$ and $i \in [t]$, verifying that $\mathrm{(iv)}$ of the lemma holds with probability $1-o(1)$. Thus (i)--(v) hold with high probability, so that the required sets $S_i$, $i\in [t]$, can be found.
\end{proof}

\section{Proof of \Cref{thm:main}} 
\label{sec: main proof}

In this section, we will prove \Cref{thm:main}. By \Cref{cor:findexpander}, it suffices to prove the following. 

\begin{theorem}
\label{edge disjoint cycles in expander}
There exists some $C$ such that for every $k$ and sufficiently large $n$ (in terms of $k$), the following holds. If $G$ is a bipartite $n$-vertex $18$-almost-regular $(\eps,s)$-expander with average degree $d(G) \ge (\log n)^C$, where
$\eps= 2^{-5}$ and $s\geq d(G)/(\log n)^2$, then $G$ contains $k$ edge-disjoint cycles with the same vertex set.
\end{theorem} 

Therefore, in the rest of this section, we prove Theorem~\ref{edge disjoint cycles in expander}. 
Throughout the proof, we will assume $C$ to be a sufficiently large constant and $n$ to be sufficiently large with respect to $k$ to support our arguments. 
Let $G$ be a bipartite $n$-vertex $18$-almost-regular $(\eps,s)$-expander with average degree $d \geq (\log n)^C$, where $\eps= 2^{-5}$ and $s\geq d/(\log n)^2$. 
Fix a bipartition $A\cup B$ for the bipartite graph $G$. Recall that a set $S\subset V(G)$ is called balanced if $|S\cap A|=|S\cap B|$. 
The following observation will be used a few times. 
\begin{enumerate}[label=(O),leftmargin=*]
    \item \label{item:balanced paths} For every $u\in A$ and $v\in B$, the internal vertex set of every $u-v$ path in $G$ is balanced,

    for every $u\in A$ and $v\in A$, every $u-v$ path in $G$ has one more internal vertex in $B$ than in $A$, and

    for every $u\in B$ and $v\in B$, every $u-v$ path in $G$ has one more internal vertex in $A$ than in $B$.
\end{enumerate}

Throughout the following proof we will note in the left margin the correspondence with the steps \ref{overview:expander}--\ref{ov:robust} in \Cref{sec:detailedoutline}, where there are some small differences to the slightly simplified outline. Note that in passing to Theorem~\ref{edge disjoint cycles in expander} we have already carried out \ref{overview:expander}.

\subsection{Finding connecting sets with good degree control}
\label{proofsubsec:nicesets}

\hspace{-1.3 cm}{\ref{overview:regsub},\ref{overview:nicesets}}

\vspace{-0.65cm}

Let $m_1:= 2n/d^{1/100}$, $m_2:= 5n/d^{1/100}$, $m_3:= 1000n/d^{1/100}$, $m_4:= n/d^{1/1000}$ and $m_5:= n/d^{1/10000}$. We now use \Cref{lem:findsets} with these values of $m_1,\ldots,m_5$ to find a subgraph $G'$ of $G$ and pairwise disjoint vertex sets $R_1, R_2, X, U_1, U_2\subseteq V(G')$ (corresponding to the sets $S_1,\dots,S_5$ in \Cref{lem:findsets}) such that the following properties are satisfied. 
Let $D:= d^{1-1/10} (\log n)^{10} $, $D_1:= d^{1-1/100} (\log n)^{10} $, and $D_2:= d^{1-1/1000} (\log n)^{10} $. 
\begin{enumerate}[label=(A\arabic*),leftmargin=*]
    \item $R_1$, $R_2$, $X$, $U_1$, and $U_2$ are balanced (with respect to the bipartition of $G$).

    \item $|R_1|=2m_1$, $|R_2|=2m_2$, $|X|=2m_3$, $|U_1|=2m_4$, and $|U_2|=2m_5$.
    
    \item \label{item:3rd} $R_1$ is $\frac{m_1^9 d}{n^9 (\log n)^{76} }$-connecting and thus also $D$-connecting since $D \le \frac{m_1^9 d}{n^9 (\log n)^{76} }$,
    
    $U_1$ is $\frac{m_4^9 d}{n^9 (\log n)^{76} }$-connecting and thus also $D_1$-connecting since $D_1\le \frac{m_4^9 d}{n^9 (\log n)^{76} }$, and 
    
    $U_2$ is $\frac{m_5^9 d}{n^9 (\log n)^{76} }$-connecting and thus also $D_2$-connecting since $D_2\le \frac{m_5^9 d}{n^9 (\log n)^{76} }$. 

    \item  \label{item:U is regular} $G'[U_1\cup U_2]$ is $\left(d'\pm 2(d')^{2/3}\right)$-nearly-regular for some $d'\ge \frac{d^{1-1/10000}}{C}$. 
    
    (To see this, note that $\mathrm{(iv)}$ of \Cref{lem:findsets} implies the existence of $d_i$ for $i\in [2]$ such that every $v\in V(G')$ satisfies $d_i- d_i^{2/3}\leq d_{G'}(v,U_i)\leq d_i+d_i^{2/3}$, where $d_2 \ge \frac{d m_5}{C n}= \frac{d^{1-1/10000}}{C}$, and we use $d_1^{2/3}+d_2^{2/3}\leq 2(d_1+d_2)^{2/3}$.)

    \item \label{item:5th} For every $v\in V(G)$, we have $d_{G}(v,R_1)\leq 2C d^{1-1/100}$, $d_{G}(v,R_2)\leq 5C d^{1-1/100}$, $d_{G}(v,X)\leq 1000C d^{1-1/100}$, and $d_{G}(v,U_1)\leq C d^{1-1/1000}$.
\end{enumerate}

\subsection{Constructing absorbers}
\label{proofsec:absorbers}

\hspace{-0.9 cm}{\ref{overview:RMBG}}

\vspace{-0.65cm}

We now use \Cref{lem:bipartite graph absorbing absorber} with $m = n/d^{1/100}$ to find a set $K\subseteq ((R_1\cup R_2)\cap A) \times ((R_1\cup R_2)\cap B)$ consisting of an odd number of pairs such that no vertex in $R_1\cup R_2$ appears in more than 102 pairs of $K$, and the following holds. 
\begin{enumerate}[label=(B),leftmargin=*]
\item \label{item:matching in R} For every balanced set $R'_1\subseteq R_1$ with $|R'_1|\le m_1 = |R_1|/2$, there is a set $K'\subseteq K$ of pairs such that $K'$ is a perfect matching in $K[(R_1\cup R_2)\setminus R'_1]$. (In other words, every vertex in $(R_1\cup R_2)\setminus R'_1$ appears in exactly one of the pairs in $K'$, and the pairs in $K'$ contain no other vertices.)
\end{enumerate}

Our goal now is to construct $k$ absorbers for every pair in $K$ (in the sense of Definition~\ref{def:absorber}). As we will see near the end of the proof, \ref{item:matching in R} allows us to show that for any balanced set $R'_1 \subseteq R_1$ of size at most $|R_1|/2$, the vertices in $(R_1 \cup R_2) \setminus R'_1$ can be paired up to be `absorbed' by our absorbers (where this pairing is given by $K'$). 

To that end, let $K \coloneqq \{p_1,\ldots,p_s\}$ where $s$ is an odd number. Note that $s \le 102 \cdot |(R_1\cup R_2)\cap A| = 102 (m_1+m_2) < 2m_3 = |X|$. Recall that $X$ is balanced. By removing vertices from $X$ if necessary, we may assume that $|X|=s+1$ and that $X$ is a balanced set. Then, by \ref{item:5th}, we still have $d_{G}(v,X)\leq 1000C d^{1-1/100}$ for every $v\in V(G)$. Enumerate the vertices in $X$ by $x_1,\ldots,x_{s+1}$ so that, for every odd $j\in [s]$, we have $x_j\in A$ and $x_{j+1}\in B$.

\hspace{-0.9 cm}{\ref{overview:absorbers}}

\vspace{-0.65cm}

Let $R= R_1\cup R_2$. Now \ref{item:5th} implies that, for every $v \in U_1$, $d_G(v, R) \le d_G(v, R_1) + d_G(v, R_2) \le 7 C d^{1 - 1/100} \le \frac{D_1}{200 k (\log n)^6}$, and, for every $v \in U_2$, $d_G(v, U_1 \cup X) \le C d^{1 - 1/1000} + 1000Cd^{1-1/100} \le \frac{D_2}{2k}$. Hence, by \ref{item:3rd}, the conditions of \Cref{finding absorber chain} are satisfied, and we obtain the following statement guaranteeing the desired absorbers.

\begin{enumerate}[label=(C),leftmargin=*] 
\item \label{item:absorber} For all $i\in [k]$ and $j\in [s]$ there is an absorber $\Abs_j^i$ in $G$ for the pair $p_j$ with endpoints $x_j$ and $x_{j+1}$ such that, writing $S_j^i$ for the interior of $\Abs_j^i$, the sets $S_j^i$ (for all $i\in [k]$ and $j\in [s]$) are pairwise disjoint subsets of $U_1\cup U_2$ of size at most $(\log n)^{12}$. 
\end{enumerate}

Let $U_{\textrm{abs}}$ be the set of vertices of $U_1\cup U_2$ used in the absorbers, i.e., $U_{\textrm{abs}} = \cup_{i\in [k], j\in [s]} S_j^i$. For every $j\in [s]$, since the vertices $x_j$ and $x_{j+1}$ are from opposite sides of the partition $A\cup B$, by \ref{item:balanced paths}, the set $S_j^i$ is balanced for every $i\in [k]$. Thus, the set $U_{\textrm{abs}}$ is balanced.
Let $U_{\textrm{unused}}=U_1\cup U_2\setminus U_{\textrm{abs}}$. Since $U_1$, $U_2$, and $U_{\textrm{abs}}$ are all balanced, the set $U_{\textrm{unused}}$ is also balanced.

\hspace{-0.9 cm}{\ref{balancing set}}

\vspace{-0.65cm}

Let $W$ be a random balanced subset of $U_{\textrm{unused}}$ of size $(1-\frac{1}{k})|U_{\textrm{unused}}|$, which is the union of a uniformly random subset of $U_{\textrm{unused}}\cap A$ of size $(1-\frac{1}{k})|U_{\textrm{unused}}\cap A|$ and a uniformly random subset of $U_{\textrm{unused}}\cap B$ of size $(1-\frac{1}{k})|U_{\textrm{unused}}\cap B|$ (here, we assume that $2k$ divides $|U_{\textrm{unused}}|$ to avoid unnecessary technicalities). 
At this stage, we point out that the vertex set of our desired $k$ (edge-disjoint) cycles will be $W\cup U_\textrm{abs}\cup X\cup R_1\cup R_2$.

\hspace{-0.9 cm}{\ref{absorbing path}}
\vspace{-0.65cm}

By \ref{item:absorber}, for every pair $p_j \in K$, there are $k$ absorbers with endpoints $x_j$ and $x_{j+1}$, namely $\Abs_j^i$, $i \in [k]$. Independently for every $j\in [s]$, generate a permutation $\sigma_j : [k] \rightarrow [k]$ uniformly at random. Now, for each pair $p_j \in K$, we will assign the absorber $\Abs_j^{\sigma_j(i)}$  to the $i$-th cycle.
Let $U_{\textrm{abs}}^i = \cup_{j\in [s]} S_j^{\sigma_j(i)}$ i.e., $U_{\textrm{abs}}^i$ is the subset of $U_{\textrm{abs}}$ consisting of those vertices which are used in the absorbers assigned to the $i$-th cycle. Then, observe that for any $i \in [k]$, the union of the absorbers assigned to the $i$-th cycle (i.e., $\Abs_j^{\sigma_j(i)}$ for $j \in [s]$) contains an $x_1$-$x_{s+1}$ path with vertex set $U_{\textrm{abs}}^i \cup X$ such that any pair in $K$ can be `absorbed' into this path (this is made precise at the end of the proof). 
Moreover, note that, for each $i\in [k]$, $U_{\textrm{abs}}^i$ is balanced. 

\subsection{Constructing {\boldmath $k$} edge-disjoint linear forests}
\label{subsec:linearforestsproof}

\hspace{-0.9 cm}{\ref{definingVi}}

\vspace{-0.65cm}

For each $i\in [k]$, let $V^i=W\cup (U_{\textrm{abs}}\setminus U_{\textrm{abs}}^i)$. Note that $V^i \subseteq U_1 \cup U_2$ (and recall that $G'[U_1\cup U_2]$ is $\left(d'\pm 2(d')^{2/3}\right)$-nearly-regular by \ref{item:U is regular}).
Let $d''=(1-\frac{1}{k})d'$. 
We will next show that the following holds with probability $1 - o(1)$. 
\begin{enumerate}[label=(D),leftmargin=*] 
\item \label{item:union of blobs are regular} For every $i\in [k]$, the graph $G'[V^i]$ is $\left(d''\pm 3(d'')^{2/3}\right)$-nearly-regular. 
\end{enumerate}
In order to prove this, fix $i\in [k]$ and consider a vertex $v\in U_1\cup U_2$. Without loss of generality, assume that $v\in (U_1\cup U_2)\cap A$ (the other case is identical).
Since \ref{item:U is regular} and \ref{item:absorber} hold, and $(\log n)^{12}\leq (d')^{1/8}$, we can apply \Cref{lem:random subgraph} for the vertex $v$ with $G=G'[U_1\cup U_2]$, $V=(U_1\cup U_2) \cap B$, $V'=U_{\textrm{unused}}\cap B$, $V^\iota_j = S^\iota_j\cap B$ for $\iota \in [k]$ and $j \in [s]$, $d=d'$, $2(d')^{2/3}$ playing the role of $d'$, $t=k$, $\tau=k-1$,  $U=W\cap B$,  $S_j=[k]\setminus \{\sigma_j(i)\}$, and $Z = (W\cup (U_{\textrm{abs}}\setminus U_{\textrm{abs}}^i)) \cap B =  V^i \cap B$. 
Then, with probability $1-n^{-\omega(1)}$, we have $d_{G'}(v, V^i) = d_{G'}(v, V^i \cap B) = d_{G'}(v, Z) = d''\pm 3(d'')^{2/3}$. 
Thus, a simple union bound over all $i\in [k]$ and $v\in U_1\cup U_2$ yields \ref{item:union of blobs are regular}.

For each $i\in [k]$, since each of $W$, $U_{\textrm{abs}}$, and $U_{\textrm{abs}}^i$ is balanced, the set $V^i$ is also balanced. 
Let $t=d^{1/5}$ and let $V^i_1,\dots,V^i_t$ be a uniformly random partition of $V^i$ into $t$ balanced sets of equal size. (Here, we assume that $2t$ divides $|V^i|$ to avoid unnecessary technicalities.) 
Note that $V^i_1\cap A,\dots,V^i_t\cap A$ is a uniformly random partition of $V^i\cap A$ into $t$ sets of equal size, and similarly $V^i_1\cap B,\dots,V^i_t\cap B$ is a uniformly random partition of $V^i\cap B$ into sets of equal size. 
Thus, since \ref{item:union of blobs are regular} holds, by a simple application of \Cref{lem:simple random subgraph} and a union bound, the following holds with probability $1 - o(1)$. For each $i\in [k]$ and each $j\in [t-1]$, the bipartite graph $G'[V^i_j,V^i_{j+1}]$ is $\left(d''/t\pm 4(d'')^{2/3}\right)$-nearly-regular. 
This, together with the inequality $6k\le (d'')^{2/3}$, implies that the following property holds with probability $1 - o(1)$. 

\begin{enumerate}[label=(E),leftmargin=*] 
\item \label{item:blobs are regular} Every graph $G''$ obtained from $G'$ by deleting at most $6k$ edges incident to every vertex has that $G''[V^i_j,V^i_{j+1}]$ is $\left(d''/t\pm 5(d'')^{2/3}\right)$-nearly-regular for each $i\in [k]$ and $j\in [t-1]$. 
\end{enumerate}
We will next show that the following statement holds with probability $1 - o(1)$. 
\begin{enumerate}[label=(F),leftmargin=*] 
\item \label{item:max degree to parts} For each $i\in [k]$, $j\in [t]$, and $v\in V(G)$, we have $d_G(v,V^i_j)\leq 25d^{4/5}$. 
\end{enumerate}
To that end, fix $i\in [k]$, $j\in [t]$, and consider a vertex $v\in A$. (Vertices in $B$ can be dealt with in a similar fashion.) 
We will apply \Cref{lem:simple random subgraph}(a) with $V=V^i\cap B$, $U=V_j^i\cap B$, and with $18d$, $0$ playing the roles of $d$, $d'$ respectively. 
By this application, since $d_G(v)\leq 18d$ and $t=d^{1/5}$, with probability $1 - n^{-\omega(1)}$ we have $d_G(v,V^i_j)=d_G(v,V^i_j\cap B)\leq 25d^{4/5}$. 
Thus, a union bound over all $i\in [k]$, $j\in [t]$ and $v\in V(G)$ yields \ref{item:max degree to parts} with probability $1 - o(1)$, as desired.

Since the events \ref{item:blobs are regular} and \ref{item:max degree to parts} hold with probability $1 - o(1)$, from now on we assume that they occur simultaneously. 

\hspace{-0.9 cm}{\ref{find matchings}}

\vspace{-0.65cm}

Our next goal is to show that we can find, for each $i\in [k]$ and $j\in [t-1]$, a matching $M^i_j$ in $G'[V^i_j,V^i_{j+1}]$ such that the following two properties hold. (Observe that for every $i\in [k]$,  the union of the matchings $M^i_j$ for $j\in [t-1]$ is a linear forest.)
\begin{enumerate}[label=(G\arabic*),leftmargin=*] 
\item \label{item:disjoint matching} For $i\in [k]$ and $j\in [t-1]$, the matchings $M^i_j$ are pairwise edge-disjoint and also edge-disjoint from all the absorbers given by \ref{item:absorber}.
\item \label{item:well behaved leftover} For every $i\in [k]$, if we let $Y^i$ be the set containing all vertices $y\in V^i_1\cup V^i_t$ that are not used in $M^i_1$ or $M^i_{t-1}$ and all vertices $y\in V^i_j$ with $2\le j\le t-1$ which are not used in both $M^i_{j-1}$ and $M^i_j$, then $|Y^i|\le 10 |V^i| d^{-1/10}\log n$ and, for every $v\in V(G)$, $|Y^i\cap N_G(v)|\leq 250d^{9/10}\log n$. 
\end{enumerate}

We will prove the existence of such matchings by induction on $i$. Assume the existence of such matchings $M^{\iota}_j$ for all $\iota < i$. To construct the matchings $M^i_1,\ldots,M^i_{t-1}$, consider the graph $G''$ that is obtained from $G'$ by deleting all the edges in $\bigcup_{\iota \in [i-1], j\in [t-1]} M^{\iota}_j$ and all the edges in $G'[U_1\cup U_2]$ that belong to one of the absorbers given by \ref{item:absorber}. Then, note that the graph $G''$ is obtained from $G'$ by deleting at most $6k$ edges incident to every vertex. (Indeed, at most $2k$ edges incident to every vertex belong to the matchings $\bigcup_{\iota \in [i-1], j\in [t-1]} M^{\iota}_j$, and at most $4$ edges incident to every vertex in $U_1 \cup U_2$ belong to the absorbers.)
Now we aim to apply \Cref{lem:matchings with little leftover} with $G$ replaced by $G''$, $V_j$ replaced by $V^i_j$ for $j \in [t]$, and $\{N_G(v): v\in V(G)\}$ playing the role of $\{T_1,\dots,T_n\}$. To that end, note that since \ref{item:blobs are regular} holds, we can apply the lemma with $\delta=d''/t-5(d'')^{2/3}$ and $\Delta=d''/t+5(d'')^{2/3}$, so that we have $1-\delta/\Delta\leq 10(d'')^{2/3}/(d''/t)= 10t/(d'')^{1/3}\leq d^{-1/9}$, where in the last inequality we used that $d'' =(1-\frac{1}{k})d'\geq d^{1-1/1000}$ holds by \ref{item:U is regular} and $t = d^{1/5}$. Hence, the condition $1-\delta/\Delta\leq t^{-1/2}\log n$ required by the lemma is satisfied. 
Therefore, since \ref{item:max degree to parts} holds, we can apply \Cref{lem:matchings with little leftover} with $r= 25d^{4/5}$ to obtain matchings $M^i_j$ in $G''[V^i_j,V^i_{j+1}]$ for $j\in [t-1]$ such that \ref{item:well behaved leftover} holds for $i$. 
By the construction of $G''$, the matchings $M^i_1,\ldots,M^i_{t-1}$ are edge-disjoint from $\bigcup_{\iota \in [i-1], j\in [t-1]} M^{\iota}_j$ and all the absorbers given by \ref{item:absorber}.
This finishes the proof of the existence of the matchings $M^i_j$ for $i\in [k]$ and $j\in [t-1]$ satisfying \ref{item:disjoint matching} and \ref{item:well behaved leftover}.

\subsection{Extending the linear forests to {\boldmath $k$} edge-disjoint cycles on the same vertex set}

\hspace{-0.9 cm}{\ref{ov:connectingpaths}}

\vspace{-0.65cm}

For $i\in [k]$ and $j\in [t-1]$, we let $\overrightarrow{Y}_j^i$ denote the vertices in $V^i_j$ that are not matched in $M^i_j$, and let $\overleftarrow{Y}_{j+1}^i$ denote the vertices in $V^i_{j+1}$ that are not matched in $M^i_j$. 
Since for every $i\in [k]$, the sets $V^i_j$ for $j\in [t-1]$ are balanced sets of equal size, we have $|\overrightarrow{Y}_j^i \cap A| = |\overleftarrow{Y}_{j+1}^i \cap B|$ and $|\overrightarrow{Y}_j^i \cap B| = |\overleftarrow{Y}_{j+1}^i \cap A|$. 
We can thus arbitrarily match the vertices in $\overrightarrow{Y}_j^i \cap A$ with the vertices in $\overleftarrow{Y}_{j+1}^i \cap B$, and the vertices in $\overrightarrow{Y}_j^i \cap B$ with the vertices in $\overleftarrow{Y}_{j+1}^i \cap A$ to obtain a collection of pairs $\mathcal{K}^i_j\subseteq \overrightarrow{Y}_j^i \times \overleftarrow{Y}_{j+1}^i$ for every $i\in [k]$ and $j\in [t-1]$. 
For every $i\in [k]$, let $\mathcal{K}^i = \cup_{j\in [t-1]} \mathcal{K}^i_j$. 
Consider the auxiliary graph $H^i$ on the vertex set $V^i$ whose edges are the pairs in $\mathcal{K}^i$. 
Then, the graph that is the union of $H^i$ and $\cup_{j\in [t-1]} M^i_j$ is a linear forest $F^i$ on $V^i$ where every path has one of its endpoints in $V^i_1$ and the other endpoint in $V^i_t$. 
For $i\in [k]$, enumerate the vertices in $V^i_1$ as $u^i_1,\ldots,u^i_\ell$ and the vertices in $V^i_t$ as $v^i_1,\ldots,v^i_\ell$ such that for every $j\in [\ell]$, there is a path in $F^i$ with endpoints $u^i_j$ and $v^i_j$. 
Let $\mathcal{K}^i_* = \{(u^i_1,x_1), (v^i_\ell,x_{s+1})\} \cup \{(u^i_{j+1}, v^i_j) : j\in [\ell-1]\}$.

Next, we will show the following.

\begin{enumerate}[label=(H\arabic*),leftmargin=*] 
\item \label{item:build long path} For every $i \in [k]$ and every pair $(u,v)$ 
in $\mathcal{K}^i\cup \mathcal{K}^i_*$, there is a path $Q^i_{uv}$ in $G$ of length at most $(\log n)^6$ with internal vertices in $R_1$ and with endpoints $u$ and $v$. Moreover, the paths $Q^i_{uv}$ (for $i \in [k]$ and $(u,v) \in \mathcal{K}^i\cup \mathcal{K}^i_*$) are pairwise internally vertex-disjoint, and  each of these paths is edge-disjoint from all the absorbers given by \ref{item:absorber}. (See Figure~\ref{fig:constructingcycle} for an illustration.)

\hspace{-1.8 cm}{\ref{path containing matchings}}

\vspace{-0.65cm}

Notice that, for each $i\in [k]$, the paths $Q^i_{uv}$ for $(u,v) \in \mathcal{K}^i\cup \mathcal{K}^i_*$ together with the edges in the matchings $M^i_1, \ldots, M^i_{t-1}$ give us a path, $P^i$ say, with endpoints $x_1$ and $x_{s+1}$. Moreover, the paths $P^i$ for $i \in [k]$ are pairwise edge-disjoint and the set of internal vertices of each path $P^i$ is $V^i\cup R^i_1$, where $R^i_1\subseteq R_1$ consists of the set of internal vertices in the paths $Q^i_{uv}$ for $(u,v) \in \mathcal{K}^i\cup \mathcal{K}^i_*$.
 
\item \label{item: at most half vertices used in R1} For every $i\in [k]$, $R^i_1$ is a balanced set containing at most $m_1$ vertices. 
\end{enumerate}

With a slight abuse of notation, let $\cup_{i\in [k]} (\mathcal{K}^i\cup \mathcal{K}^i_*)$
denote the multiset obtained by taking the union of the sets $\mathcal{K}^i\cup \mathcal{K}^i_*$ for $i\in [k]$. To prove \ref{item:build long path}, we first use (from \ref{item:3rd}) that $R_1$ is $D$-connecting to construct $5$ paths of length at most $(\log n)^6$ between each pair of the multiset $\cup_{i\in [k]} (\mathcal{K}^i\cup \mathcal{K}^i_*)$, with internal vertices in $R_1$ such that all these paths are pairwise internally vertex-disjoint. To achieve this, consider the multiset $S$ obtained by replacing each pair in $\cup_{i\in [k]} (\mathcal{K}^i\cup \mathcal{K}^i_*)$ with $5$ copies of it. Note that any vertex of $G$ appears at most $2\cdot k \cdot 5=10k$ times in this multiset $S$. 
Also, observe that the pairs in $\cup_{i\in [k]}(\mathcal{K}^i\cup \mathcal{K}^i_*)$ do not contain any vertex outside of $\{x_1,x_{s+1}\}\cup \left(\cup_{i\in [k]} (V^i_1 \cup V^i_t \cup Y^i)\right)$. Thus, by using \ref{item:max degree to parts} and \ref{item:well behaved leftover}, every vertex in $G$ has at most $10 k \cdot (k \cdot 250d^{9/10}\log n+k \cdot 2 \cdot 25d^{4/5}+2) \le 10k^2 \cdot 251 d^{9/10}\log n \le d^{1-1/10} (\log n)^{10} = D$ neighbours in $S$. 
So, by \Cref{def:nice}, we have $5$ paths of length at most $(\log n)^6$ between each pair  $(u,v) \in \cup_{i\in [k]} (\mathcal{K}^i\cup \mathcal{K}^i_*)$ with internal vertices in $R_1$ such that all these paths are pairwise internally vertex-disjoint. 
On the other hand, for every pair $(u,v) \in \cup_{i\in [k]} (\mathcal{K}^i\cup \mathcal{K}^i_*)$, among these $5$ pairwise internally vertex-disjoint paths between the pair $(u,v)$, at most $4$ paths are not edge-disjoint from the absorbers given by \ref{item:absorber}. Indeed, all these paths have internal vertices in $R_1$, so the only edges of the absorbers that can appear in these paths are between $U_{\textrm{abs}}$ and $R_1$, but each vertex in $U_{\textrm{abs}}$ is incident to at most $2$ edges (with the other endpoint in $R_1$) that are used in the absorbers. This shows that the desired path $Q_{uv}^i$ exists between every pair $(u,v) \in \mathcal{K}^i\cup \mathcal{K}^i_*$ for every $i \in [k]$, verifying \ref{item:build long path}. 

To prove \ref{item: at most half vertices used in R1}, we fix $i\in [k]$. 
Since \ref{item:well behaved leftover} holds, the number of pairs in $\mathcal{K}^i\cup \mathcal{K}^i_*$ is at most $|Y^i| + |V^i_1| + |V^i_t| \le 10nd^{-1/10}\log n + 2nd^{-1/5} \le 11nd^{-1/10}\log n$. Therefore, since each of the paths $Q_{uv}^i$ given by \ref{item:build long path} (between the pairs $(u,v) \in \mathcal{K}^i\cup \mathcal{K}^i_*$ for $i \in [k]$) have length at most $(\log n)^6$, we have $|R^i_1| \le 11nd^{-1/10}\log n \cdot (\log n)^6 \le m_1$. 
By construction, every pair $(u,v) \in \mathcal{K}^i$ contains one vertex from $A$ and the other from $B$. This ensures that the internal vertex sets of the paths $Q_{uv}^i$ given by \ref{item:build long path}  between these pairs are balanced sets. 
Moreover, note that the set $\{x_1,x_{s+1}\} \cup V^i_1 \cup V^i_t$ is balanced (recall that $x_1\in A$ and $x_{s+1}\in B$), and every vertex in this set appears in exactly one pair of $\mathcal{K}^i_*$. 
Thus, by \ref{item:balanced paths}, the union of the internal vertex sets of the paths between the pairs in $\mathcal{K}^i_*$ is a balanced set. Thus, it follows that $R^i_1$ is a balanced set. 
This finishes the verification of \ref{item: at most half vertices used in R1}. 

\hspace{-0.9 cm}{\ref{ov:robust}}

\vspace{-0.65cm}

Finally, for every $i\in [k]$, our aim is to extend the path $P^i$ given by \ref{item:build long path} to a cycle $C^i$ with the vertex set $W\cup U_{\text{abs}}\cup X\cup R_1\cup R_2$ such that the cycles $C^i$ (for $i \in [k]$) are pairwise edge-disjoint. 
To that end, by \ref{item:build long path}, for every $i\in [k]$, it suffices to find a path $P^i_*$ with endpoints $x_1$ and $x_{s+1}$ and with $V(P^i_*) = U^i_{\textrm{abs}} \cup X \cup (R_1\setminus R^i_1) \cup R_2$ such that the paths $P^i_*$ (for $i \in [k]$) are
pairwise edge-disjoint and do not contain any edge outside of the absorbers given by \ref{item:absorber}.
(Indeed, then for every $i \in [k]$, by concatenating $P^i$ and $P^i_*$ we obtain a cycle $C^i$ with the vertex set $W\cup U_{\text{abs}}\cup X\cup R_1\cup R_2$ such that the cycles $C^i$ for $i \in [k]$ are pairwise edge-disjoint, as required.)
To find the paths $P^i_*$, for every $i\in [k]$, let $K^i\subseteq K$ be the set of pairs of vertices obtained by applying \ref{item:matching in R} with $R^i_1$ playing the role of $R'_1$ (this application is possible since \ref{item: at most half vertices used in R1} holds). Then, every vertex of $(R_1\setminus R^i_1) \cup R_2$ appears in exactly one of the pairs in $K^i$ and the pairs in $K^i$ contain no other vertices. 
We now fix an arbitrary $i\in [k]$. 
Consider the set $J^i\subseteq [s]$ such that $K^i = \{p_j : j\in J^i\}$. 
For every $j\in J^i$, consider the path $L_j$ with endpoints $x_j$ and $x_{j+1}$ and with internal vertex set $S_j^{\sigma_j(i)} \cup p_j$ which is contained in the absorber $\Abs_j^{\sigma_j(i)}$. 
For every $j\in [s]\setminus J^i$, consider the path $L_j$ with endpoints $x_j$ and $x_{j+1}$ and with internal vertex set $S_j^{\sigma_j(i)}$ which is contained in the absorber $\Abs_j^{\sigma_j(i)}$. 
Now, by concatenating the paths $L_1,\ldots,L_s$, we obtain the desired path $P^i_*$ with endpoints $x_1$ and $x_{s+1}$ and with $V(P^i_*) = U^i_{\textrm{abs}} \cup X \cup (R_1\setminus R^i_1) \cup R_2$. 
Moreover, since the sets $S^i_j$ for $i\in [k]$ and $j\in [s]$ are pairwise disjoint and every edge in $\Abs_j^i$ is incident to a vertex in $S_j^i$, the paths $P^i_*$ are pairwise edge-disjoint, as required. 
This shows the existence of $k$ edge-disjoint cycles on the common vertex set $W\cup U_{\text{abs}}\cup X\cup R_1\cup R_2$, completing the proof of \Cref{edge disjoint cycles in expander}.

\section{Connecting pairs of vertices by vertex-disjoint paths through a random vertex set} 
\label{sec: sublinear expander}

In this section we prove \Cref{lem:randomsetisnice}. As mentioned in the introduction, we will use an adaption of arguments by Buci\'c and Montgomery \cite{bucic2022towards} (which uses, in part, some ideas from Tomon~\cite{MR4695962}). Before we turn to the details of the proof, it is worth highlighting how our result and argument differ from those in \cite{bucic2022towards}. The result in \cite{bucic2022towards} closest to our \Cref{lem:randomsetisnice} is \cite[Theorem 16]{bucic2022towards}, which roughly states the following: if $G$ is a good robust (sublinear) expander, and $V$ is a $1/3$-random subset of $V(G)$, then, given any collection of pairs in $V(G)$ in which each vertex appears only a small number of times, we can connect all the pairs using short \emph{edge-disjoint} paths through $V$. Note that our \Cref{lem:randomsetisnice} differs from this result in a number of ways. First, crucially, we are looking to find \emph{vertex-disjoint}, rather than \emph{edge-disjoint} paths which connect the pairs. This requires us to impose a stronger condition on the pairs that we seek to connect, since it is not possible to connect more than $|V|$ pairs through $V$ using internally vertex-disjoint paths. More precisely, in \cite{bucic2022towards} the condition on the pairs is that each vertex should appear only a small number of times in the pairs, while we require that each vertex in the random set $V$ has few neighbours among the vertices we seek to connect. Finally, we also need a generalisation from $1/3$-random sets to $p$-random sets. The argument in \cite{bucic2022towards} did not crucially rely on the set being $1/3$-random, so it is straightforward to adapt their argument to the more general $p$-random setting. On the other hand, in order to find vertex-disjoint, rather than edge-disjoint paths, we need to modify the argument in \cite{bucic2022towards} nontrivially.

A crucial result in the proof of Theorem 16 in \cite{bucic2022towards} states that, for an expander $G$, if we take a large random subset $V$ in $V(G)$, then with high probability it is true that for every $U\subset V(G)$ and not too large $F\subset E(G)$, more than half of the vertices in $V$ can be reached from $U$ using short paths inside $V$, avoiding all edges in $F$. We will need a similar, but quantitatively stronger result here (see \Cref{lem:reachable} below). In order to prove such a strengthening, the main new ingredient is \Cref{lem:partitionedgesintoexpanders}, which states that any $(\eps,s)$-expander $G$ can be edge-decomposed into almost $s$ weaker expanders. This strengthens the corresponding result in \cite{bucic2022towards} (Lemma 15) where it was shown that one can decompose the edges of $G$ into almost $\sqrt{s}$ weaker expanders. A more substantially new result that we will need and which has no analogue in \cite{bucic2022towards} is \Cref{lem:connectafterexpand}.
A detailed explanation of why we need this tweak compared to the argument in \cite{bucic2022towards} is provided before the statement of \Cref{lem:connectafterexpand}. In general, where possible we use similar constants to \cite{bucic2022towards} for ease of comparison.

\subsection{Proof of \Cref{lem:randomsetisnice}}

Given $U, V \subseteq V (G)$, the ball of radius $i$ around $U$ within $V$, denoted by $B_G^i(U, V )$, is the set of vertices in $V$ which can be reached by a path of length at most $i$ starting from a vertex in $U$ which has all of its internal vertices in $V$. The starting vertex in $U$ is not required to be in $V$ itself. We do, however, only consider reachable vertices within $V$, so that $B_G^i(U, V) \subseteq V$.

It is convenient to use the following definition.

\begin{defn}[$\lambda$-reachable set]\label{defn:reachable}
    Let $G$ be an $n$-vertex graph. We say that a set $V\subset V(G)$ is $\lambda$-reachable if, for every $U\subset V(G)$ and every $F\subset E(G)$ with $|F|\leq \lambda |U|$,
    $$|B_{G-F}^{(\log n)^4} (U,V)|>\frac{|V|}{2}.$$
\end{defn}

Most of the work in this section will go into proving the following lemma.

\begin{lemma} \label{lem:reachable}

Let $0 < p < 1$. Suppose that $G$ is an $n$-vertex $(\eps, s)$-expander with $2^{-90} < \eps < 1$ and $s\ge p^{-4}(\log n)^{30} $. Let $V$ be a $p$-random subset of $V(G)$.
Then, with probability $1-o(1)$, $V$ is $\lambda$-reachable for $\lambda=\frac{p^8s}{(\log n)^{60}}$.
\end{lemma}

In the rest of this subsection, we will complete the proof of \Cref{lem:randomsetisnice} assuming \Cref{lem:reachable}. Subsections~\ref{subsec:expanderproperties}--\ref{subsec:improved robustness} are then devoted to proving \Cref{lem:reachable}.

We will need the following, which is \cite[Proposition 8]{bucic2022towards}.

\begin{prop} \label{prop:connect one of t}
    Let $1\leq \ell,t\leq n$. Let $G$ be an $n$-vertex graph and let $V\subset V(G)$ with $|V|\geq 4t-2$ be such that, for every $U\subset V(G)$ with size $|U|=t$, we have $|B^{\ell}_G(U,V)|>\frac{|V|}{2}$. Let $z_1,\dots,z_{2t-1},w_1,\dots,w_{2t-1}$ be distinct vertices of $G$. Then, for some $j\in [2t-1]$, there is a $z_j-w_j$ path in $G$ with internal vertices in $V$ and with length at most $4\ell \log n$.
\end{prop}

We will also use the following form of the Aharoni-Haxell hypergraph matching theorem (see Corollary~1.2 in \cite{aharoni2000hall}).

\begin{theorem}\label{thm:hyperhall}
Let $r\in \N$, and let $H_1,\ldots,H_r$ be a collection of hypergraphs with at most $\ell$ vertices in each edge. Suppose that, for each $I\subseteq [r]$, there is a matching in $\bigcup_{i\in I}H_i$ containing more than $\ell(|I|-1)$ edges. Then, there is an injective function $f:[r]\to \bigcup_{i\in [r]}E(H_i)$ such that $f(i)\in E(H_i)$ for each $i\in [r]$ and $\{f(i):i\in [r]\}$ is a matching of $r$ edges.
\end{theorem}

Let us now informally discuss how we will complete the proof of \Cref{lem:randomsetisnice}, and how the argument differs from the one in \cite{bucic2022towards}. Let $V$ be a $p$-random subset of $V(G)$ and let $(x_1,y_1),\dots,(x_r,y_r)$ be the set of pairs that we seek to connect using short internally vertex-disjoint paths through $V$. Similarly to \cite{bucic2022towards}, the existence of such paths will be proven by applying \Cref{thm:hyperhall}, with $H_i$ defined to be the hypergraph whose edges correspond to the short paths through $V$ connecting $x_i$ and $y_i$. (Since here we are looking for vertex-disjoint paths, the edges of $H_i$ are the \emph{vertex sets}, rather than the \emph{edge sets}, of the paths connecting $x_i$ and $y_i$.) To verify that the assumptions in \Cref{thm:hyperhall} are satisfied, we need to prove that for each $I\subseteq [r]$ there are many vertex-disjoint short paths through $V$ that connect $x_j$ to $y_j$ for some $j\in I$. If we were looking for edge-disjoint paths, we could apply \Cref{prop:connect one of t} repeatedly: even after finding a large collection $P_1,P_2,\dots,P_q$ of paths connecting some pair $(x_j,y_j)$ (with $j\in I$), we can find a new one that is edge-disjoint from all the previous ones by applying \Cref{prop:connect one of t} with $\{z_1,\dots,z_{2t-1}\}=\{x_j:j\in I\}$, $\{w_1,\dots,w_{2t-1}\}=\{y_j:j\in I\}$ and $G-F$ in place of $G$, where $F$ is the set of edges used by the paths $P_1,\dots,P_q$. It is crucial here that if $q$ is not already large, then $F$ is quite small, so (using \Cref{lem:reachable}) the condition $|B^{\ell}_{G-F}(U,V)|>\frac{|V|}{2}$ in \Cref{prop:connect one of t} is indeed satisfied for every set $U$ of size roughly $|I|$. Unfortunately, this approach (which is the one used in \cite{bucic2022towards}) does not extend to the case where we are looking for vertex-disjoint paths, because in order to guarantee that the obtained paths are vertex-disjoint, we need to let $F$ contain every edge in $G$ that is adjacent to a vertex in one of the paths $P_1,\dots,P_q$. This makes $F$ too large for the condition $|B^{\ell}_{G-F}(U,V)|>\frac{|V|}{2}$ to hold for every set $U$ of size roughly $|I|$.

To deal with this issue, we use the following lemma, which allows us to construct the connecting paths in two stages: first, using the bounds on the degrees of vertices from $V$ to the set of pairs that we seek to connect, we prove that the set of pairs `expand' via a large star matching into a random subset $W\subset V$ of size roughly $|V|/3$. This allows us to use the argument sketched above if we instead apply Proposition~\ref{prop:connect one of t} to a larger set of pairs chosen from the leaves of these stars.

\begin{lemma} \label{lem:connectafterexpand}
    Let $n$ be sufficiently large and let $G$ be a graph with maximum degree $\Delta$ and let $V\subset V(G)$ be a $\lambda$-reachable subset.
    Let $x_1,\dots, x_{r}, y_1, \dots, y_{r}$ be a sequence of (not necessarily distinct) vertices outside of $V$. Let $W\subset V(G)\setminus (V\cup \{x_1,\dots , x_{r}, y_1, \dots , y_{r}\})$ be such that $|W|\leq |V|$ and assume that there exists some $\delta_0$ such that each $x_i$ and $y_i$ sends at least $\delta_0$ edges to $W$ and each vertex in $W$ has at most $\frac{\lambda \delta_0}{2\Delta(\log n)^{12}}$ neighbours in the multiset $\{x_1,\dots,x_{r},y_1,\dots,y_{r}\}$.
    
    Then, there is a collection of internally vertex-disjoint $x_i-y_i$ paths (one for each $i\in [r]$) of length at most $(\log n)^6$ with internal vertices in $W\cup V$.
\end{lemma}

\begin{proof}
    We first prove the following claim using the degree conditions.
    \begin{claim*}
        There exist pairwise disjoint sets $X_1,\dots,X_{r},Y_1,\dots,Y_{r}\subset W$ of size at least $\lambda^{-1}2\Delta(\log n)^{12}$ such that $X_i\subset N_G(x_i)$ and $Y_i\subset N_G(y_i)$ for all $i\in [r]$.
    \end{claim*}

    \begin{proof}[ of Claim]
        Let $q=\lambda^{-1} 2 \Delta (\log n)^{12}$. Note that, as $V$ is $\lambda$-reachable, we have $\lambda\leq \Delta$, so $q\geq 2(\log n)^{12}$. Let us assume for simplicity that $q$ is an integer. Define a bipartite graph $H$ with parts $A:=\{x_i^j: i \in [r],j\in [q]\}\cup \{y_i^j: i \in [r],j\in [q]\}$ and $W$, where the vertices $x_i^j$ and $y_i^j$ are all distinct, and in which there is an edge between $x_i^j$ and some $w\in W$ if and only if $x_i w\in E(G)$, and similarly there is an edge between $y_i^j$ and some $w\in W$ if and only if $y_i w\in E(G)$. By the assumption on the degrees in $G$, each $u\in A$ has degree at least $\delta_0$ in $H$, and each $w\in W$ has degree at most $q\cdot \frac{\lambda \delta_0}{2\Delta (\log n)^{12}}=\delta_0$ in $H$. Hence, by Hall's theorem, there exists an injective map $f:A\rightarrow W$ such that, for each $u\in A$, there is an edge in $H$ between $u$ and $f(u)$. For each $i\in [r]$, let $X_i=f(\{x_i^j: j\in [q]\})$ and $Y_i=f(\{y_i^j:j\in [q]\})$. It is straightforward to verify that these sets satisfy the conditions in the claim.
    \end{proof}
    
    For each $i\in [r]$, let $H_i$ be the hypergraph on vertex set $V$ in which a hyperedge corresponds to the set of internal vertices of a path of length at most $(\log n)^6$ between $X_i$ and $Y_i$ with all internal vertices in $V$. Note that in order to prove the lemma, it suffices to find a matching of $r$ edges in which the $i$-th edge belongs to $H_i$ for each $i\in [r]$.

    Using \Cref{thm:hyperhall}, it suffices to verify that, for each $I\subseteq [r]$, there is a matching in $\cup_{i\in I} H_i$ containing more than $(|I|-1)(\log n)^6$ edges. For this, let $I\subseteq [r]$ and let $M_I$ be a maximal matching in $\cup_{i\in I} H_i$. We need to show that $|M_I|>(|I|-1)(\log n)^6$. Assume for a contradiction that $|M_I|\leq (|I|-1)(\log n)^6$. Let $S$ be the set of vertices in $G$ used by the edges in $M_I$. Note that $|S|\leq |M_I|(\log n)^6 \leq |I|(\log n)^{12}$. Let $F$ be the set of edges in $G$ which are incident to at least one vertex in $S$. Then $|F|\leq |S|\Delta\leq |I|\Delta (\log n)^{12}$. Let $t=\lambda^{-1}|I|\Delta (\log n)^{12}$. Now, if $U\subset V(G)$ and $|U|=t$, then $|F|\leq \lambda |U|$, so as $V$ is $\lambda$-reachable, we have
    $$|B_{G-F}^{(\log n)^4} (U,V)|>\frac{|V|}{2}.$$
    Note also that, by the claim above, we have $|V|\geq |W|\geq 2r\cdot \lambda^{-1}2\Delta (\log n)^{12}\geq 4t$. Hence, we can apply \Cref{prop:connect one of t} with $G-F$ in place of $G$ and with $\ell=(\log n)^4$ to conclude that if $z_1,\dots,z_{2t-1},w_1,\dots,w_{2t-1}$ are distinct vertices in $G$, then for some $j\in [2t-1]$ there is a $z_j$-$w_j$ path in $G-F$ with internal vertices in $V$ and with length at most $4(\log n)^5\leq (\log n)^6$. However, since $|X_i|,|Y_i|\geq \lambda^{-1}2\Delta (\log n)^{12} $, we have $|\cup_{i\in I} X_i|\geq |I|\lambda^{-1}2\Delta (\log n)^{12} =2t$ and $|\cup_{i\in I} Y_i|\geq 2t$, so there is a path in $G-F$ with internal vertices in $V$ and with length at most $(\log n)^6$ which connects some element of $X_i$ to some element of $Y_i$, for some $i\in I$. Since this path does not use the edges in $F$, it does not have any internal vertex which is in $S$. Hence, the internal vertices of this path are disjoint from the vertices used by $M_I$, contradicting the maximality of $M_I$.
\end{proof}

We are now ready to prove \Cref{lem:randomsetisnice} (assuming \Cref{lem:reachable}, which will be proved in the rest of this section).

\begin{proof}[ of \Cref{lem:randomsetisnice}]
    Let $W$ be a random subset of $V$ obtained by including each vertex of $V$ independently at random with probability $1/3$. Let $V'=V\setminus W$. Note that $V'$ is a $(2p/3)$-random subset of $V(G)$, so by \Cref{lem:reachable}, $V'$ is $\lambda$-reachable with probability $1-o(1)$, where $\lambda=\frac{(2p/3)^8 s}{(\log n)^{60} }$. Let $v$ be an arbitrary vertex in $G$. Since $W$ is a $(p/3)$-random subset of $V(G)$, the expected number of neighbours of $v$ in $W$ is at least $\delta p/3$. Hence, by the lower bound on $p$ and by the Chernoff bound, the probability that $v$ has fewer than $\delta p/6$ neighbours in $W$ is $o(1/n)$. Hence, with probability $1-o(1)$, every vertex in $G$ has at least $\delta p/6$ neighbours in $W$. Also, the probability that $|V'|\geq |W|$ is $1-o(1)$.
    
    We now show that, if $V'$ is $\lambda$-reachable, each $v\in V(G)$ has at least $\delta p/6$ neighbours in $W$ and $|V'|\geq |W|$, then $V$ is $D$-connecting for $D = \frac{p^9 s \delta}{\Delta (\log n)^{73}}$.
    
    Let $x_1,\dots,x_r,y_1,\dots,y_r$ be a sequence of vertices outside of $V$ and suppose that every $v\in V$ has at most $D$ neighbours in the multiset $\{x_1,\dots,x_r,y_1,\dots,y_r\}$.
    Let $\delta_0=\delta p/6$. Now note that each $x_i$ and $y_i$ sends at least $\delta_0$ edges to $W$ and each vertex in $W$ has at most $D\leq \frac{\lambda \delta_0}{2\Delta (\log n)^{12} }$ neighbours in the multiset $\{x_1,\dots,x_r,y_1,\dots,y_r\}$. Hence, by \Cref{lem:connectafterexpand} (applied with $V'$ in place of $V$), there is a set of internally vertex-disjoint $x_i-y_i$ paths (one for each $i\in [r]$) of length at most $(\log n)^6$ with internal vertices in $W\cup V'=V$. Thus, $V$ is indeed $D$-connecting, completing the proof. 
\end{proof}

It remains to prove \Cref{lem:reachable}.

\subsection{Properties of the expander}
\label{subsec:expanderproperties}

For a graph $G$, it will be convenient to define the `robust neighbourhood' of a set $U \subseteq V(G)$ for any parameter $d$ as
$N_{G,d}(U) \coloneqq \{v \in V(G) \setminus U : |N_G(v) \cap U| \ge d \}$,
i.e., the set of vertices in $G$, outside of the set $U$, which have degree at least $d$ in $U$. We will use Proposition 12 from \cite{bucic2022towards}.

\begin{prop} \label{prop:expansion-a-b}
Let $G$ be an $n$-vertex $(\eps,s)$-expander, $U\subseteq V(G), |U|\le \frac23 n$ and $F$ a set of at most $s|U|/2$ edges. Then, for any $0<d \le s$, either
\[
\text{\textbf{\emph{a)}} }\;\;|N_{G-F}(U)| \ge \frac{s|U|}{2d},\;\;\;\text{ or }\;\;\;\text{\textbf{\emph{b)}}}\;\; |N_{G-F,d}(U)|\ge \frac{\eps|U|}{(\log n)^{2}}.
\]
\end{prop}

The following lemma is a generalisation of Proposition 13 from \cite{bucic2022towards} which shows that more structure can be found in both outcomes of the above proposition. The proof is essentially the same.

\begin{lemma} \label{lem:structuredichotomy}
    There is an $n_0$ such that the following holds whenever $n\geq n_0$, $0< \eps< 1$, $r\geq (\log n)^2$, $t\geq (\log n)^2$ and $s\geq 20rt$. Let $G$ be an $n$-vertex $(\eps, s)$-expander, let $U \subset V (G)$ have size $|U| \leq 2n/3$ and let $F$ be a set of at most $s|U|/4$ edges. Then, in $G - F$ we can find either
    \begin{enumerate}[label=\alph*)]
        \item $\frac{|U|}{10r}$ vertex-disjoint stars, each with $t$ leaves, centre in $U$ and all leaves in $V(G)\setminus U$, or \label{prop:large stars}
        \item a bipartite subgraph $H$ with vertex classes $U$ and $X\subset V(G)\setminus U$ such that \label{prop:robust nhood}
        \begin{itemize}
            \item $|X|\geq \frac{\eps |U|}{2(\log n)^2}$ and
            \item every vertex in $X$ has degree at least $r$ in $H$ and every vertex in $U$ has degree at most $2t$ in $H$.
        \end{itemize}
    \end{enumerate}
    
\end{lemma}

\begin{proof}
Take a maximal collection of vertex-disjoint stars in $G-F$ with $t$ leaves, centre in $U$ and leaves outside of $U$. Let $C \subseteq U$ be the set of centres of these stars and $L\subseteq V(G)\setminus U$ be the set consisting of all their leaves. Assuming that \ref{prop:large stars} does not hold, we have $|C| \le \frac{|U|}{10r}$ and $|L| \le |C|\cdot t \le \frac{|U|}{10r}\cdot t$, and, by the maximality, there is no vertex in $U\setminus C$ with at least $t$ neighbours in $G-F$ in $V(G)\setminus (U\cup L)$. Thus,
\begin{equation}\label{eqn:NGW}
|N_{G-F}(U\setminus C)|\leq |C|+|L|+|U\setminus C|\cdot t\leq  \frac{|U|}{10 r}+|C|\cdot t +|U\setminus C|\cdot t< 2|U|\cdot t.
\end{equation}

We now construct the set $X\subseteq V(G)\setminus U$ and the bipartite subgraph $H$ through the following process, starting with $X_0=\emptyset$ and setting $H_0$ to be the graph with vertex set $U\cup X_0$ and no edges. Let $k=|V(G)\setminus U|$ and label the vertices of $V(G)\setminus U$ arbitrarily as $v_1,\ldots,v_k$. For each $i\geq 1$, if possible, pick a star $S_i$ in $G-F$ with centre $v_i$ and $r$ leaves in $U$ such that these leaves in $U$ have degree at most $2t$ in the graph $H_{i-1}\cup S_i$, and let $H_i=H_{i-1}\cup S_i$ and $X_i=X_{i-1} \cup \{v_i\}$; otherwise set $H_i=H_{i-1}$ and $X_{i}=X_{i-1}$. Finally, let $H=H_k$ and $X=X_k=V(H_k)\setminus U$. We will show that \ref{prop:robust nhood} holds for this choice of $H$ with bipartition $(U,X)$.

Firstly, observe that every vertex in $U$ has degree at most $2t$ in $H_i$ for each $i\in [k]$ by construction, and that every vertex $v_i$ in $X$ has degree exactly $r$ in $H$, so the second condition in \ref{prop:robust nhood} holds. Thus, we only need to show that $|X| \ge \frac{\eps|U|}{2(\log n)^2 }$ holds, which will follow as no vertex in $U\setminus C$ has $t$ neighbours in $G-F$ in $X\setminus L$.

Indeed, let $U'$ be the set of vertices in $U \setminus C$ with degree exactly $2t$ in $H$. As each vertex in $U\setminus C$ has fewer than $t$ neighbours in $G-F$ in $X\setminus L$, the vertices in $U'$ must have at least $t$ neighbours in $H$ in $X\cap L$. As each vertex in $X\cap L$ has $r$ neighbours in $H$, we have
\[
|U'|\leq \frac{r|X \cap L|}{t}\le \frac{r}{t}\cdot |L| \leq \frac{r}{t}\cdot \frac{|U|\cdot t}{10r}= \frac{|U|}{10}.
\]

Let $B=C \cup U'$, so that
\[
|B|\leq \frac{|U|}{10r}+\frac{|U|}{10}\le \frac{|U|}{2},
\]
and, thus, $|U\setminus B|\geq \frac{|U|}2$.

Then, by \Cref{prop:expansion-a-b} applied to $U\setminus B$ and $F$ with $d=r$, using that $|F| \le s|U|/4 \le s|U \setminus B|/2$, we have either $|N_{G-F}(U\setminus B)|\geq \frac{s|U\setminus B|}{2r}$ or $|N_{G-F,r}(U\setminus B)|\geq \frac{\eps |U\setminus B|}{(\log n)^2}$. As $|U\setminus C|\leq |U|$ and
$$\frac{s|U\setminus B|}{2r}\geq \frac{s|U|}{4r}\geq 5t|U|,$$ $|N_{G-F}(U\setminus B)|\geq \frac{s|U\setminus B|}{2r}$ would contradict \eqref{eqn:NGW}, and therefore we must have that $|N_{G-F,r}(U\setminus B)|\geq \frac{\eps |U\setminus B|}{(\log n)^2}$. 
Every vertex $v$ in $N_{G-F,r}(U\setminus B)$ has at least $r$ neighbours in $G-F$ in $U\setminus B$ which must all have degree strictly less than $2t$ in $H$ (as they are not in $B=U' \cup C$, and using the definition of $U'$). This implies $v \in X$, since we could add it together with some $r$ of these neighbours to the graph $H$. Hence, we must have $N_{G-F,r}(U\setminus B)\subseteq X$, and 
\[
|X|\geq |N_{G-F,r}(U\setminus B)|\geq \frac{\eps |U\setminus B|}{(\log n)^2 }\geq  \frac{\eps|U|}{2(\log n)^2},
\]
as required.
\end{proof}

The next lemma is a straightforward variant of Proposition 18 from \cite{bucic2022towards}.

\begin{lemma}
\label{lem:findingwell-expandingset}
    Let $n\geq 2$, $0<\eps<1$ and $s\geq \lambda\geq 1$. Let $G$ be an $n$-vertex $(\eps,s)$-expander and let $U\subset V(G)$ have size $|U|\leq 2n/3$.
    Then, there is a set $U'\subset U$ with $|N_G(U')|\geq |U'|\lambda$ and $|U'|\geq \frac{\eps |U|}{3\lambda (\log n)^2}$.
\end{lemma}
\begin{proof}
Let $U'\subseteq U$ be maximal subject to $|N_{G}(U')|\geq |U'|\lambda$, noting that such $U'$ exists as $U'=\emptyset$ satisfies these conditions. Suppose that $U\neq U'$, for otherwise $U$ satisfies the conditions itself. Then $|N_{G}(U')|< (|U'|+1)\lambda$ or we could add an arbitrary vertex to $U'$ and contradict the maximality of $U'$. Similarly we know that, for every vertex $v\in U\setminus U'$, $v$ has at most $\lambda$ neighbours outside of $U'\cup N_{G}(U')$, for otherwise $U'\cup \{v\}$ contradicts the maximality of $U'$. Let $F$ be the set of edges between $U\setminus U'$ and $V(G)\setminus (U'\cup N_{G}(U'))$, so that $|F|\leq |U\setminus U'|\lambda\leq s|U|$. Thus, we have by the definition of expansion that 
\[
\frac{\eps|U|}{(\log n)^2} \le |N_{G-F}(U)| \le |N_{G}(U')| \le (|U'|+1)\lambda,
\]
so that $|U'|\geq \frac{\eps |U|}{3 \lambda (\log n)^{2}}$, as required.
\end{proof}

\subsection{Edge decomposition of an expander into weaker expanders}
\label{subsec:edge-decompositionexpanders}

In this subsection, we prove the following lemma which states that every expander can be edge-partitioned into many weaker expanders. Such a result was proved in \cite{bucic2022towards} (see their Lemma 15), but for our application it is crucial to have a quantitatively stronger statement.

\begin{lemma}\label{lem:partitionedgesintoexpanders} Let $n$ and $s$ be sufficiently large, $k\in \N$ and $0<\eps < 1$.
Suppose that $G$ is an $n$-vertex $(\eps,s)$-expander and $\frac{\eps s}{k} \geq 10^5 (\log n)^3$. Then, there are edge-disjoint graphs $G_1,\ldots,G_k$ such that $E(G)=\bigcup_{i\in [k]}E(G_i)$ and, for each $i\in [k]$, $G_i$ is an $\left(\frac{\eps}4,\frac{\eps s}{10^4k(\log n)^2}\right)$-expander  with vertex set $V(G)$.
\end{lemma}

In order to obtain this stronger statement, we will use the following simple proposition instead of using Proposition~\ref{prop:expansion-a-b} as is done in \cite{bucic2022towards}.

\begin{prop}\label{prop:expansion-red-blue}
Let $n$ and $s$ be sufficiently large and let $0<\eps < 1$. Let $G$ be an $n$-vertex $(\eps,s)$-expander, and let $U\subseteq V(G)$ with $|U|\le \frac23 n$. Then, there is some $0\leq i\leq \log s$ such that
\begin{equation}\label{eqn:somei}
|N_{G,2^i}(U)|\geq \frac{\eps s|U|}{2^i(\log n)^2}.
\end{equation}
\end{prop}
\begin{proof} Suppose for a contradiction that \eqref{eqn:somei} does not hold for any $0\leq i\leq \log s$.
Let $U'$ be the set of vertices in $V(G)\setminus U$ with at most $s$ neighbours in $U$. Let $F$ be the set of edges between $U$ and $U'$ in $G$, and note that 
\[
|F|\leq \sum_{i=0}^{\log s} |\{u\in U':2^i\leq d_{G}(u,U)\leq 2^{i+1}\}|\cdot 2^{i+1}\leq \sum_{i=0}^{\log s}|N_{G,2^i}(U)|\cdot 
2^{i+1}\leq \sum_{i=0}^{\log s}2^{i+1}\cdot \frac{\eps s|U|}{2^i(\log n)^2}\leq s|U|,
\]
where we have used that $s\leq n$ follows from the definition of expansion.
Then, since $G$ is an $(\eps,s)$-expander, we have $|N_{G,s}(U)| \ge |N_{G-F}(U)|\geq \eps|U|/(\log n)^2$, so \eqref{eqn:somei} holds for $i=\log s$, which is a contradiction.
\end{proof}

The following result states that if we randomly sample the edges of an expander, then with high probability we get a (weaker) expander. We will use it to prove \Cref{lem:partitionedgesintoexpanders}.

\begin{lemma}\label{lem:edgesampleexpander}
    Let $n$ and $s$ be sufficiently large and $0<p,\eps<1$.
Suppose that $G$ is an $n$-vertex $(\eps,s)$-expander and $\eps p s \geq 10^5 (\log n)^3$. Let $H$ be a random subgraph of $G$ with vertex set $V(G)$, which contains every edge independently with probability $p$. Let $s'=\frac{\eps p s}{10^4(\log n)^2}$. Then, the probability that $H$ is not an $(\frac{\eps}{4},s')$-expander is less than $2/n$.
\end{lemma}

\begin{proof}
Let $U\subseteq V(G)$ and $u=|U|\leq \frac{2n}{3}$. We will show that the probability that $U$ does not satisfy the $(\eps/4,s')$-expansion condition is at most $e^{-2u\log n}$. Applying \Cref{prop:expansion-red-blue},
there is some $0\leq i\leq \log s$ such that, letting $d=2^i$, $|N_{G,d}(U)|\geq \frac{\eps s u}{d\cdot (\log n)^2}$. 

First consider the case \textbf{a)} where $pd\leq 100$. Then, for each $v\in N_{G,d}(U)$, the probability that $v\in N_{H}(U)$ is at least $1-(1-p)^{d}\geq \frac{pd}{10^3}$. Note that,  therefore, $|N_{H}(U)|$ is dominated by $\bin(|N_{G,d}(U)|, \frac{pd}{10^3})$. Hence $\E[|N_{H}(U)|]\geq \frac{\eps p s u}{10^3 (\log n)^2} \ge 100 u \log n$. Thus, by the Chernoff bound, with probability at least $1 - e^{-2u\log n}$, we have $|N_{H}(U)|\geq \frac{\eps p s u}{2\cdot 10^3 (\log n)^2} \geq 2s'u$, which implies that for any set $F\subset E(H)$ with $|F|\leq s'u$, we have $|N_{H-F}(U)| \geq s'u \geq 10u\log n \ge \frac{\eps u}{4 (\log n)^2}$. 

Next, consider the case \textbf{b)} where $pd> 100$. Then, for each $v\in N_{G,d}(U)$, the probability that $v\notin N_{H,pd/2}(U)$ is at most 
$q \coloneqq \P(\bin(d, p) < \frac{pd}{2})$. By the Chernoff bound, $q \le e^{-pd/12}$.
Writing $t=\frac{\eps s u}{2d(\log n)^2}$, as $|N_{G,d}(U)|\geq \frac{\eps su}{d(\log n)^2}=2t$, we have
\begin{align}\label{eqn:NHU}
\P\left(|N_{H,pd/2}(U)|<t\right)&\le \binom{\ceil{2t}}{\ceil{2t}-\floor{t}}\cdot q^{\ceil{2t}-\floor{t}}  \le 2^{\ceil{2t}}\cdot q^{\ceil{2t}-\floor{t}}\le (4q)^t \le e^{\frac{-pdt}{24}}= e^{-\frac{\eps p su}{48 (\log n)^2}}\notag\\&\le e^{-2u\log n}.
\end{align}
Note that, if $|N_{H,pd/2}(U)|\geq \frac{\eps s u}{2d(\log n)^2}$, then, for any $F\subset E(H)$ with $|F|\leq s'u$, we have 
\[
|N_{H-F}(U)|\geq |N_{H,pd/2}(U)|-\frac{|F|}{pd/2}\geq \frac{\eps s u}{2d(\log n)^2}-\frac{2s'u}{pd}\geq \frac{\eps s u}{4d(\log n)^2}\geq  \frac{\eps u}{4(\log n)^2}.
\]

Therefore, whichever of \textbf{a)} or \textbf{b)} holds, the probability that $U$ does not satisfy the $(\eps/4,s')$-expansion condition is at most $e^{-2u\log n}$. Hence, the probability that $H$ is not an $(\eps/4,s')$-expander is at most
$$ \sum_{u=1}^{2n/3}\binom{n}{u}e^{-2u\log n} \le \sum_{u=1}^{2n/3} n^u\cdot n^{-2u} \le n^{-1}+\sum_{u=2}^{2n/3} n^{-2}< \frac{2}{n},$$
as required.
\end{proof}

\begin{proof}[ of \Cref{lem:partitionedgesintoexpanders}]
Assign every edge of $G$ to one of the graphs $G_1,\ldots, G_k$ uniformly and independently at random, so that every $G_i$ is a random subgraph of $G$ containing each edge of $G$ with probability $1/k$. For any given $i\in [k]$, by \Cref{lem:edgesampleexpander}, the probability that $G_i$ is not an 
$\left(\frac{\eps}4,\frac{\eps s}{10^4k(\log n)^2}\right)$-expander is strictly less than $1/k$ (since $n\geq s\geq 2k)$.  Thus, by a union bound, the probability that 
$G_i$ is an $\left(\frac{\eps}4,\frac{\eps s}{10^4k(\log n)^2}\right)$-expander for every $i \in [k]$ is positive, so some decomposition as required by the lemma must exist.
\end{proof}

\subsection{Expansion of well-expanding sets into a random vertex set}

In order to prove \Cref{lem:reachable}, we need to show that, for an expander $G$, if we take a large random subset $V$ in $V(G)$, then with high probability it is true that, for every $U\subset V(G)$ and not too large $F\subset E(G)$, more than half of the vertices in $V$ can be reached from $U$ by short paths inside $V$ which avoid all the edges in $F$. We prove this in three steps: in this subsection we deal only with `well-expanding' sets $U$ and very small $F$, in the next subsection we extend this to arbitrary $U$ but still only very small $F$, and finally, in Subsection~\ref{subsec:improved robustness}, we deal with all $U$ and much larger $F$, completing the proof of \Cref{lem:reachable}.
Similar results were proved in \cite{bucic2022towards}, but we need quantitatively stronger versions here.

\begin{lemma}
\label{lem:wellexpandingsetscanreach}
    Let $n$ be sufficiently large, let $0 < p < 1$, and suppose that $G$ is an $n$-vertex $(\eps, s)$-expander with $2^{-100} < \eps < 1$ and $s \geq 20p^{-3}(\log n)^{13}$. Let $U \subset V(G)$ satisfy $|N_G(U)| \geq |U|p^{-4}(\log n)^{24}$ and let $F \subset E(G)$ satisfy $|F| \leq |U|$. Let $V$ be a $p$-random subset of $V(G)$.
    Then, with probability $1-e^{-\Omega(|U|(\log n)^2)}$,
    $$|B_{G-F}^{(\log n)^4} (U,V)|>\frac{|V|}{2}.$$
\end{lemma}

\begin{proof}
Let $\ell=(\log n)^4$ and let $q\in (0,1)$ be such that $1-(1-q)^{\ell-1}(1-\frac{4p}{5})=p$, i.e., that $(1-q)^{\ell-1}=\frac{1-p}{1-\frac{4p}{5}}$, so that
\begin{equation}\label{eqn:p15}
q\ge \frac{p}{6(\log n)^4}.
\end{equation}
Independently, for each  $i\in [\ell]$, let $V_i$ be a $q$-random subset of $V(G)$ if $i\leq \ell-1$ and a $(4p/5)$-random subset of $V(G)$ if $i=\ell$. Set $V=V_1\cup\ldots \cup V_\ell$, and note that $V$ is a $p$-random subset of $V(G)$. Thus, we wish to show that, with probability at least  $1-e^{-\Omega\left({|U|}{(\log n)^{2}}\right)}$ we have $|B^{\ell}_{G-F}(U , V)|>\frac{|V|}{2}$.

For each $0\le i\leq \ell$, let $B_{i}$ be the set of vertices of $G$ which can be reached via a path in $G-F$ which starts in $U$ and has length at most $i$ and whose internal vertices (if there are any) are in $V_1\cup \dots \cup V_{i-1}$. In particular, we have $B_0=U$ and $B_1=U\cup N_{G-F}(U)$. 
Observe also that $B_0\subseteq B_1\subseteq \ldots \subseteq B_{\ell}$. We emphasise that the vertices of $B_{i}$ do not themselves have to be inside $V_1\cup \ldots \cup V_{i-1}$, only the internal vertices of some path from $U$ to the vertex in $B_{i}$ are required to be inside $V_1\cup \ldots \cup V_{i-1}$. An important property of $B_{i}$ is that it is completely determined by the sets $U,V_1,\ldots, V_{i-1}$, and therefore is independent of $V_{i}$. Note also that any vertex in $N_{G-F}(B_i)$ with a neighbour in $B_i$ that gets sampled into $V_i$ belongs to $B_{i+1}$. These two observations will be the key behind why the sets $B_{i+1}$ will grow in size until they occupy most of the set $V(G)$. The lemma will then follow from
\begin{equation}\label{eqn:overkill4}
B_\ell\cap V_\ell\subseteq B^{\ell}_{G-F}(U, V).
\end{equation}

We now show that indeed, for each $1\le i\le \ell-1$, unless $B_{i}$ is already very large, $B_{i+1}$ is likely to be somewhat larger than $B_i$.

\begin{claim*} For each $1\le i\le \ell-1$, with probability $1-e^{-\Omega\left({|U|}{(\log n)^{2}}\right)}$, either $|B_i| \ge \frac23 n$, or  $$|B_{i+1}\setminus B_i| \ge \frac{\eps|B_i|}{2^5(\log n)^{2} }.$$
\end{claim*}
\begin{proof} For each $v\in N_{G-F}(B_i)$, $v$ is in $B_{i+1}$ if at least one of its neighbours in $G-F$ in $B_i$ gets sampled into $V_{i}$. That is,
\begin{equation}\label{eqn:overkill1}
\{v\in N_{G-F}(B_i):(N_{G-F}(v)\cap B_i)\cap V_i\neq\emptyset\} \subseteq B_{i+1}\setminus B_i. 
\end{equation}
We will show that, for any set $W\subseteq V(G)$ with $|W|\leq \frac23 n$ and $B_1\subseteq W$ 
\begin{equation}\label{eqn:overkill2}
\P\left(|\{v\in N_{G-F}(W):(N_{G-F}(v)\cap W)\cap V_i\neq\emptyset\}|\geq \frac{\eps|W|}{2^5(\log n)^{2}}\right)\ge 1-e^{-\Omega\left({p^4|B_1|}/{(\log n)^{22} }\right)}.
\end{equation}
Given~\eqref{eqn:overkill2}, we will have that for all $1\le i \le \ell -1$,
\begin{align*}
\P\left(|B_i| \ge \frac23 n  \: \text{ or  } \: |B_{i+1}\setminus B_i| \ge \frac{\eps|B_i|}{2^5(\log n)^{2}}\right) &\overset{\textcolor{white}{\eqref{eqn:overkill1}}}{\geq} \P\left(|B_{i+1}\setminus B_i| \ge \frac{\eps|B_i|}{2^5(\log n)^{2}}\: \Big| \: |B_i|\leq \frac{2}{3}n\right)
\\
&\overset{\eqref{eqn:overkill1}}{\geq} 
\P\left(|\{v\in N_{G-F}(B_i):(N_{G-F}(v)\cap B_i)\cap V_i\neq\emptyset\}|\ge \frac{\eps|B_i|}{2^5(\log n)^{2}}\: \Big| \: |B_i|\leq \frac{2}{3}n\right)\\
&\overset{\eqref{eqn:overkill2}}{\ge} 
1-e^{-\Omega\left({p^4|B_1|}/{(\log n)^{22} }\right)}\ge 1-e^{-\Omega\left(|U| (\log n)^{2} \right)},
\end{align*}
where in the last inequality we used that $|B_1|\geq |N_{G-F}(U)|\ge |U|p^{-4}(\log n)^{24}-|F|\geq \frac12|U|p^{-4}(\log n)^{24}$.

Let then $W\subseteq V(G)$ satisfy $|W|\leq \frac23n$ and $B_1\subseteq W$. As $|W|\leq \frac23 n$ and $|F|\le |U|\leq |B_1|\leq |W|\le {s|W|}/4$, we can apply \Cref{lem:structuredichotomy} with $W$ in place of $U$, $r=p^{-1}(\log n)^4$ and $t=p^{-2}(\log n)^9$ (note that the lemma applies since $s \ge 20 p^{-3} (\log n)^{13} = 20 rt$). Hence, one of the two cases \ref{prop:large stars} or \ref{prop:robust nhood} from \Cref{lem:structuredichotomy} holds; we will show that \eqref{eqn:overkill2} holds in either case.

a) Suppose that $G-F$ contains $\frac{|W|}{10p^{-1}(\log n)^{4}}$ vertex-disjoint stars, each with $p^{-2}(\log n)^9$ leaves, centre in $W$ and all leaves in $N_{G-F}(W)$. Let $C\subseteq W$ be the set of centres of such a collection of stars, and note that
\begin{equation}\label{eqn:overkill3}
|\{v\in N_{G-F}(W):(N_{G-F}(v)\cap W)\cap V_i\neq \emptyset\}|\geq |C\cap V_i|\cdot p^{-2}(\log n)^9.
\end{equation}
By the Chernoff bound and \eqref{eqn:p15}, and as $p^2/(\log n)^{14}\leq 1$ and $|W|\geq |B_1|$, with probability  at least 
$1-e^{-q|C|/12}= 1-e^{-\Omega\left({|W|p^2}/{(\log n)^{8}}\right)} \geq 1-e^{-\Omega\left({p^4|B_1|}/{(\log n)^{22}}\right)}$, we have $|C\cap V_i|\geq \frac{q|C|}2\geq \frac{|W|p^2}{120(\log n)^{8}}$. This, combined with \eqref{eqn:overkill3}, implies \eqref{eqn:overkill2}.

b) Suppose instead that there is a bipartite subgraph $H\subseteq G-F$ with vertex classes $W$ and $X\subseteq V(G)\setminus W$ such that
    \begin{itemize}
        \item $|X| \ge \frac{\eps|W|}{2(\log n)^2}$ and
        \item every vertex in $X$ has degree at least $r=p^{-1}(\log n)^{4} $ in $H$ and every vertex in $W$ has degree at most $\Delta:=2t=2p^{-2}(\log n)^{9}$ in~$H$.
    \end{itemize}  
    For each $v\in X$, the probability that $v$ has no neighbours in $H$ in $V_i$ is at most
    \[
(1-q)^r= (1-q)^{p^{-1}(\log n)^4}\le e^{-qp^{-1}(\log n)^4}\overset{\eqref{eqn:p15}}{\leq} e^{-1/6}\leq \frac{7}{8}.
    \]
Let $Y$ be the random variable counting the number of vertices in $X$ having a neighbour in $V_{i}$ in $H$, so that $\E [Y] \ge \frac{|X|}{8}$. Observe also that $Y$ is $\Delta$-Lipschitz since for each $v\in W$ the event $\{v\in V_i\}$ affects $Y$ by at most $d_H(v)\leq \Delta$. 
Hence, by \Cref{lem:mcd} with $k=\Delta$, $t=\frac{|X|}{16}$ and $N=|W|$, we have
$$\P\left(Y < \frac{|X|}{16}\right)\le \P\left(Y < \E [Y] - \frac{|X|}{16} \right)\le2\exp\left(-\frac{2^{-9}|X|^2}{\Delta^2|W|}\right)= e^{-\Omega\left({|W|p^4}/{(\log n)^{22}}\right)}. $$
Each vertex in $X$ with a neighbour in $V_i$ in $H$ belongs to $\{v\in N_{G-F}(W):(N_{G-F}(v)\cap W)\cap V_i\neq\emptyset\}$. Hence, with probability at least $1-e^{-\Omega\left({|W|p^4}/{(\log n)^{22}}\right)}$, we have $|\{v\in N_{G-F}(W):(N_{G-F}(v)\cap W)\cap V_i\neq\emptyset\}|\geq Y\geq \frac{|X|}{16}\geq \frac{\eps|W|}{2^5(\log n)^2}$. This means that \eqref{eqn:overkill2} holds in case b) as well, completing the proof of the claim.
\end{proof}

As $B_\ell$ and $V_\ell$ are independent, by the Chernoff bound, we have that
\[
\P\left(|B_{\ell}\cap V_\ell|\le\frac{31p}{60}n \: \big| \: |B_{\ell}|\geq \frac{2}3n\right)\leq \P\left(\bin\left(\frac{2}3n,\frac{4p}5\right)\leq \frac{31p}{60}n\right)\le e^{-\Theta(np)},
\]
and, similarly, we have $\P\left(|V|\geq \frac{61p}{60}n\right)\le e^{-\Theta(np)}$.

Thus, by the claim, altogether we have that

\begin{enumerate}[label=\roman*)]
\item \label{itm1}
for each $i\in [\ell-1]$, $|B_i| \ge \frac23 n$ or  $|B_{i+1}\setminus B_i| \ge \frac{\eps|B_i|}{2^5(\log n)^{2}}$, 

 \item $|B_\ell|<\frac23 n$ or $|B_\ell\cap V_\ell|>\frac{31p}{60}n$, and

\item \label{itm3} $|V|\leq \frac{61p}{60}n$
\end{enumerate}

with probability at least
$$1-(\log n)^4\cdot e^{-\Omega\left({|U|}{(\log n)^{2}}\right)}-e^{-\Theta(np)}\ge 1-e^{-\Omega\left({|U|}{(\log n)^{2} }\right)},$$ where we used that $|U|\le \frac{np^4}{(\log n)^{24}}$ as $|N_G(U)|\geq |U|p^{-4}(\log n)^{24}$.

However, if \ref{itm1}--\ref{itm3} all hold, then, for each $i\in [\ell]$, we have
$$|B_i| \ge \min\left\{\frac23 n,\left(1+\frac{\eps}{2^5(\log n)^2 }\right)^{i}|U|\right\}\geq \min\left\{\frac23 n,\exp\left(\frac{\eps i}{2^6(\log n)^2}\right)\right\},$$
so that, setting $i=\ell=(\log n)^4$, we obtain that $|B_{\ell}|\geq \frac23 n$, and hence, by ii) and iii), that $|B_\ell\cap V_\ell|>\frac{|V|}2$.

Thus, by \eqref{eqn:overkill4}, we have that $|B^{\ell}_{G-F}(U,V)|>\frac{|V|}2$ with probability at least $1-e^{-\Omega\left({|U|}{(\log n)^{2}}\right)}$.
\end{proof}

\subsection{Expansion of all sets into a random vertex set with poor robustness}
We now use \Cref{lem:findingwell-expandingset} and \Cref{lem:wellexpandingsetscanreach} to prove the following result.   
\begin{lemma} \label{lem:reachhalfwithfewedgesblocked}
    Let $0 < p < 1$. Suppose that $G$ is an $n$-vertex $(\eps, s)$-expander with $2^{-100} < \eps < 1$ and $s \geq 20 p^{-4}(\log n)^{24}$. Let $V$ be a $p$-random subset of $V(G)$.
    Then, with probability $1-o(1/n)$, for every $U\subset V(G)$ and every $F\subset E(G)$ with $|F|\leq \frac{p^4|U|}{(\log n)^{27}}$,
    \begin{equation}
    \label{eqn:B1UV}
    |B_{G-F}^{(\log n)^4} (U,V)|>\frac{|V|}{2}.    
    \end{equation}
\end{lemma}
\begin{proof} 
Say that a set $U'\subseteq V(G)$ is \emph{well-expanding} in $G$ if 
$|N_{G}(U')| \ge |U'|p^{-4}(\log n)^{24}.$ Since $s \geq 20p^{-3}(\log n)^{13}$, given a non-empty well-expanding set $U'\subseteq V(G)$ and a set of edges $F$ of size
at most $|U'|$, \Cref{lem:wellexpandingsetscanreach} applied to $U'$ implies that 
\begin{equation}\label{eq:1}
    |B_{G-F}^{(\log n)^4}(U', V)|> \frac{|V|}2
\end{equation}
fails with probability at most $e^{-\Omega(|U'|(\log n)^2)}.$

Now, a union bound over all pairs $(U',F)$ such that $U'$ is a well-expanding set in $G$ and $F$ is a set of at most $|U'|$ edges tells us that \emph{some} such pair $(U',F)$ fails \eqref{eq:1} with probability at most 
\begin{align*}
    \sum_{(U',F)}  e^{-\Omega(|U'|(\log n)^2)} &\le 
    \sum_{u=1}^n\sum_{f=1}^{u} \binom{n}{u}\binom{n^2}{f} \cdot e^{-\Omega(u(\log n)^2)} \\
&    \le \sum_{u=1}^n u\cdot n^{3u}\cdot e^{-\Omega(u(\log n)^2)} \le \sum_{u=1}^n e^{-\Omega(u(\log n)^2)}=o(1/n).
\end{align*}

Thus, with probability $1-o(1/n)$, we can assume that \eqref{eq:1} holds for every well-expanding set $U'$ and set $F\subseteq E(G)$ with $|F|\leq |U'|$. We will now show that this implies that \eqref{eqn:B1UV} holds for all $U\subseteq V(G)$ and $F\subseteq E(G)$ with $|F|\leq \frac{p^4|U|}{(\log n)^{27}}$, completing the proof.

Let $U\subseteq V(G)$ with $|U|\leq \frac23n$ and let $F\subseteq E(G)$ satisfy the (slightly weaker) condition $|F|\leq \frac{2p^4|U|}{(\log n)^{27}}$. Then, applying \Cref{lem:findingwell-expandingset} (with $\lambda = p^{-4}(\log n)^{24} \le s$), there is a set $U'\subseteq U$ which is well-expanding for which $|U'|\geq \frac{\eps p^4 |U|}{3(\log n)^{26}}$. Noting that $|F|\leq |U'|$ (as we may assume that $n$ is large), we therefore have that
\[
|B_{G-F}^{(\log n)^4}(U,V)|\geq |B_{G-F}^{(\log n)^4}(U',V)|> \frac{|V|}{2}.
\]
Finally, consider $U\subseteq V(G)$ with $|U|> \frac23n$ and let $F\subseteq E(G)$ satisfy $|F|\leq \frac{p^4|U|}{(\log n)^{27}}$. Let $\bar{U}\subseteq U$ be an arbitrary subset with $\frac{n}2\leq |\bar{U}|\leq \frac23 n$, so that we have $|F|\leq \frac{2p^4 |\bar{U}|}{(\log n)^{27}}$, and hence, from what we have just shown,
\[
|B_{G-F}^{(\log n)^4}(U,V)|\geq |B_{G-F}^{(\log n)^4}(\bar{U},V)|> \frac{|V|}{2},
\]
as required.
\end{proof}

\subsection{Expansion of all sets into a random vertex set with improved robustness} \label{subsec:improved robustness}

We are now ready to prove \Cref{lem:reachable} in the following (equivalent) form which incorporates Definition~\ref{defn:reachable}.

\begin{lemma} \label{lem:reachhalfwithmanyedgesblocked}

Let $0 < p < 1$. Suppose that $G$ is an $n$-vertex $(\eps, s)$-expander with $2^{-90} < \eps < 1$ and $s\ge p^{-4}(\log n)^{30} $. Let $V$ be a $p$-random subset of $V(G)$.
Then, with probability $1-o(1)$, for every $U\subset V(G)$ and every $F\subset E(G)$ with $|F|\leq \frac{p^8|U|s}{(\log n)^{60}}$,
    $$|B_{G-F}^{(\log n)^4} (U,V)|>\frac{|V|}{2}.$$
\end{lemma}

\begin{proof}
    Let $k= \lfloor p^4s/(\log n)^{30}  \rfloor$. Then by the assumption on $s$, we have $k \ge 1$. Also, as our lemma concerns the asymptotic probability as $n\rightarrow \infty$, we may assume that $n$ is sufficiently large and therefore $\eps s/k \ge 10^5 (\log n)^3$. Hence, using Lemma \ref{lem:partitionedgesintoexpanders}, we can obtain edge-disjoint graphs $G_1,\dots,G_k$ such that $E(G)=\cup_{i\in [k]} E(G_i)$ and, for each $i\in [k]$, $G_i$ is an $(\frac{\eps}{4},s')$-expander with vertex set $V(G)$, where $s'=\frac{\eps s}{10^4 k(\log n)^2}$. By the choice of $k$, we have $s'\geq 20 p^{-4}(\log n)^{24}$, so by Lemma \ref{lem:reachhalfwithfewedgesblocked} and a union bound over the $k$ graphs $G_i$, with probability $1-o(1)$ we have that, for each $i\in [k]$ and every $U\subset V(G_i)$ and $F\subset E(G_i)$ with $|F|\leq \frac{p^4|U|}{(\log n)^{27}}$,
    $$|B_{G_i-F}^{(\log n)^4} (U,V)|>\frac{|V|}{2}.$$

    Now let $U\subset V(G)$ and $F\subset E(G)$ with $|F|\leq \frac{p^8|U|s}{(\log n)^{60}}$. As the graphs $G_i$, $i\in [k]$, are edge-disjoint, there exists some $i\in [k]$ such that $|F\cap E(G_i)|\leq \frac{p^8|U|s}{k(\log n)^{60}} \le \frac{p^4|U|}{(\log n)^{27}}$, and therefore
    $$|B_{G-F}^{(\log n)^4} (U,V)|\geq |B_{G_i-(F\cap E(G_i))}^{(\log n)^4} (U,V)|>\frac{|V|}{2}.$$
    This proves the lemma.
\end{proof}

\section{Concluding remarks} \label{sec:concluding}

Note that our regularisation lemma (\Cref{thm:effnearreg2}) is vacuous when the average degree $d$ is less than about $\log n$. However, using the same approach, one can also obtain meaningful analogues of \Cref{thm:effnearreg2} when $d$ is very small. Indeed, one can replace the use of the union bound by an application of the Lov\'asz local lemma in the proof of \Cref{lem:onestephighprob}, and obtain the following results. (Note that a downside is that, due to the use of the local lemma, the results no longer have the useful property that, for example, with high probability the nearly-regular subgraph contains a large random subset of $V(G)$.)

\begin{lemma}\label{lem:allverticesregularspanningnearregular}
    There exists some $C>0$ such that for each $1/d\ll \gamma\leq 1/100$ the following holds. Let $G$ be a graph with $d \leq \delta(G), \Delta(G) \leq (1+\gamma)d$. Then, for some $d'\geq (1-40\gamma)d$, $G$ contains a subgraph $G'$ with $|V(G')|\geq (1-40\gamma)|V(G)|$ and degrees between $d'$ and $d'+C\log d'$.
\end{lemma}

\begin{lemma}\label{thm:nearreg}
    There exists some $t\geq 1$ such that the following holds for all $C\geq 1$ and $C'\geq C^t$. For any $d$, if a graph has degrees between $d$ and $Cd$, then it contains a subgraph with degrees between $d'$ and $d'+C'\log d'$ for some $d'\geq d/C'$.
\end{lemma}

We remark that \Cref{lem:allverticesregularspanningnearregular} and \Cref{thm:nearreg} will be used by Montgomery, M\"uyesser, Pokrovskiy and Sudakov in their upcoming work \cite{MMPSpathdecomp}.
Moreover, in our upcoming paper \cite{CJMMregular}, using further ideas, we prove that every graph with degrees between $d$ and $C d$ contains an $r$-regular subgraph, where $r\geq d/C'$ for some $C'$ which is a polynomial in $C$. This shows that the conclusion of \Cref{thm:nearreg} can be strengthened to find a fully regular subgraph. Furthermore, this can be used to make significant progress on an old problem of R\"odl and Wysocka \cite{rodl1997note} (see \cite{CJMMregular}). \\

\textbf{Acknowledgements.} We are grateful to Noga Alon and Benny Sudakov for helpful comments and suggestions.

\bibliographystyle{abbrv}
\bibliography{bib}

\appendix
\appendixpage

\section{Finding an 18-almost-regular expander}
In this appendix we prove \Cref{cor:findexpander}. For convenience, we repeat the lemma along with the definition of our expander.
\defnsublinear*
\findexpander*

We will use the following lemma, which finds a robust expander in any graph $G$ whose average degree is very close to that of $G$, while also ensuring that the minimum degree of the expander is large.

\begin{lemma}\label{lem:findexpander}
    Let $n$ be a sufficiently large integer and $0< \eps < 2^{-3}$. Let $G$ be an $n$-vertex graph with $d(G)\geq (\log n)^2$. Then, $G$ contains a subgraph $G'$ satisfying the following properties.
    \begin{itemize}
        \item $G'$ is an $(\eps,s)$-expander for some $s\geq \frac{d(G')}{(\log |V(G')|)^2}$.
        \item $d(G')\geq d(G)\big(1 - \frac{50}{\log \log n}\big)$.
        \item $\delta(G') \ge \frac{d(G)}{3}$.
    \end{itemize}
\end{lemma}

\begin{proof}
We will define a procedure that finds the desired subgraph $G'$ in $G$. At every step, we consider a subgraph $H$ and show that either $H$ satisfies the desired properties (in which case we end the procedure and set $G':=H$) or we can find a certain subgraph $H' \subseteq H$ and continue the procedure with $H'$. Before describing the procedure, we prove the following claim.

\begin{claim}\label{claim:dichotomy}
Let $H$ be a subgraph of $G$ which is not an $(\eps,s)$-expander with $s = \frac{d(H)}{(\log |V(H)|)^2}$, and let $\lambda = \frac{3}{(\log |V(H)|)^2}$. Then there is a set $Y \subseteq V(H)$ with $|Y| < \frac{3}{4} |V(H)|$ such that $d(H[Y]) \ge (1 - \lambda)d(H)$ or there is a set $X\subseteq V(H)$ with $X \neq V(H)$ such that $d(H[X]) \ge d(H)$.
\end{claim}

\begin{proof}[ of \Cref{claim:dichotomy}]
As $H$ is not an $(\eps,s)$-expander, there are sets $U \subseteq V(H)$ and $F \subseteq E(H)$ with $1 \le |U| \le \frac{2}{3} |V(H)|$ and $|F| \le s|U|$ such that $|N_{H-F}(U)| < \frac{\eps|U|}{(\log |V(H)|)^2}$. Let $Y = U \cup N_{H-F}(U)$ and $X = V(H) \setminus U$. 
Note that 
\begin{equation}
\label{eqn:splitH}
d(H) (|U| + |X|) = 2 e(H) \le 2e(H[Y]) + 2e(H[X])+ 2 |F| \le 2e(H[Y]) + 2e(H[X])+ 2 s |U|.
\end{equation}

We claim that either $d(H[Y]) \ge (1 - \lambda)d(H)$ or $d(H[X]) \ge d(H)$. Indeed, otherwise, by \eqref{eqn:splitH}, we have 
$$d(H)(|U| +|X|) \le (1 - \lambda)d(H)|Y| + d(H)|X| + 2 s|U| \le (1 - \lambda)d(H)(|U|+ |N_{H-F}(U)|) + d(H)|X| + 2s |U|.$$

Hence, we have
$$d(H)(|U| +|X|) \le (1 - \lambda)d(H) \bigg(|U|+ \frac{\eps |U|}{(\log |V(H)|)^2} \bigg) + d(H)|X| + 2s |U|.$$
Rearranging, and dividing by $|U|$, we obtain
\begin{equation}
\label{eq:relatedands}
 d(H) \bigg( \lambda - \frac{(1 - \lambda) \eps}{(\log |V(H)|)^2}\bigg) \le 2s.   
\end{equation}

However, by the choice of $\lambda$, the left-hand side of \eqref{eq:relatedands} is more than $d(H) \left( \lambda - \frac{\eps}{(\log |V(H)|)^2}\right) = d(H) \left( \frac{3 - \eps}{(\log |V(H)|)^2}\right) > 2s$, a contradiction. Therefore, we have either $d(H[Y]) \ge (1 - \lambda)d(H)$ or $d(H[X]) \ge d(H)$. Moreover, since
$\eps < 2^{-3}$, in the former case, $|Y| = |U \cup N_{H-F}(U)| < \big(1 + \frac{\eps}{(\log |V(H)|)^2}\big) |U| \le \big(1 + \frac{\eps}{(\log |V(H)|)^2}\big)\cdot \frac{2}{3} |V(H)| < \frac{3}{4} |V(H)|$, as desired. This proves the claim.
\end{proof}

Now, we describe our procedure. Let $d_1=d(G)$. Starting with $G_1 \coloneqq G$, iteratively, for every $i\ge 1$, we do the following as long as $|V(G_i)|\ge d_1/2$ holds, where we use $\lambda = \frac{3}{(\log |V(H)|)^2}$. 
\begin{itemize}
\item If $G_i$ has a vertex $v$ with degree less than $d(G_i)/2$, we define $G_{i+1} := G_i \setminus v$ and proceed to the next step with $G_{i+1}$. (Note that, in this case, $d(G_{i+1}) \ge d(G_i)$.) Otherwise, we proceed to the next bullet point. 

\item If $G_i$ is an $(\eps,s)$-expander with $s = \frac{d(G_i)}{(\log |V(G_i)|)^2}$, then define $G'=G_i$ and stop this procedure. Otherwise, we apply \Cref{claim:dichotomy} with $H=G_i$ to get a set $X$ or $Y$ as in the claim. In the former case, we have $d(H[X]) \geq d(H)$, and so define $G_{i+1}:= H[X]$ and proceed to the next step with $G_{i+1}$. For the latter case, we have $d(H[Y]) \geq (1-\lambda)d(H)$, and so define $G_{i+1}:= H[Y]$ and proceed to the next step with $G_{i+1}$, noting that $|V(G_{i+1})| = |Y|  < \frac{3}{4} |V(G_i)|$.
\end{itemize}

Suppose that the above procedure terminates with the graph $G_t$. For every $i\in [t]$, define $n_i:=|V(G_i)|$ and $d_i:=d(G_i)$.
By the stopping condition, note that $n_{t-1}\ge d_1/2$.
Our goal now is to show that the graph $G_t$ is our desired subgraph $G'$. Note that, in the above procedure, we have $d_{i+1} < d_i$ only if $n_{i+1} < \frac{3}{4} n_i$. Moreover, by \Cref{claim:dichotomy}, in this case, $d_{i+1} \ge (1 - \frac{3}{(\log n_i)^2}) d_i$. We now wish to obtain a lower bound on $d_t$. Let $I$ be the set of $i\in [t-1]$ for which $d_{i+1}<d_i$, noting that, as $n_{i+1}<\frac{3}{4}n_i$ for each $i\in I$, for each $j$ there is at most one value of $i\in I$ with $(d_1/2)\cdot (4/3)^j\leq n_i<(d_1/2)\cdot (4/3)^{j+1}$.   
Then, we have 
\begin{equation}
\label{eq:boundingdt}
d_t\ge d_1 \cdot \prod_{i\in I} \left(1 - \frac{3}{(\log n_i)^2}\right) 
\ge d_1\cdot \prod_{i\ge 0}\left(1 - \frac{3}{(\log \left(d_1/2\cdot (4/3)^i\right))^2}\right) 
\ge d_1\cdot \bigg(1-\sum_{i\ge 0} \frac{3}{(\log \left(d_1/2 \cdot (4/3)^i\right))^2}\bigg). 
\end{equation}
Let $k$ be an integer such that $(4/3)^k \le d_1/2 \le (4/3)^{k+1}$. Then, since $n$ is large enough, $$\sum_{i\ge 0} \frac{3}{(\log \left(d_1/2 \cdot (4/3)^i\right))^2} \le \sum_{i\ge 0} \frac{3}{(\log \left((4/3)^{k+i}\right))^2} = \sum_{i\ge k} \frac{3}{i^2 (\log (4/3))^2} \le \sum_{i\ge k} \frac{20}{i^2} \le \frac{40}{k} \le \frac{100}{\log d_1}.$$
Thus, by \eqref{eq:boundingdt}, we have $$d_t \ge d_1 \left(1 - \frac{100}{\log d_1}\right) \ge d(G) \left(1 - \frac{50}{\log \log n}\right).$$

Since $d_t \ge d_1/2$, we have $n_t \ge d_1/2$. Therefore, the procedure could have only stopped because $G' \coloneqq G_t$ is an $(\eps,s)$-expander with $s = \frac{d(G_t)}{(\log |V(G_t)|)^2}$. Moreover, it has average degree $d(G_t) = d_t \ge d(G) \left(1 - \frac{50}{\log \log n}\right)$, and minimum degree at least $d(G_t)/2 \ge d(G)/3$. Thus, $G' = G_t$ is the desired subgraph. This completes the proof of the lemma.
\end{proof}

We combine the above lemma with the following result, which can be proved using an approach of Pyber~\cite{pyber1985regular}, as is done by Alon, Cohen, Dey, Griffiths, Musslick, Ozcimder, Reichman, Shinkar and Wagner~in \cite{ARS+17} and by Buci\'c, Kwan, Pokrovskiy, Sudakov, Tran and Wagner~in \cite{bucic2020nearly}.

\begin{lemma} \label{lem:almostreg} Every $n$-vertex graph with average degree $d$ contains a $6$-almost-regular subgraph with average degree at least $\frac{d}{100\log n}$.    
\end{lemma}

We are now ready to prove \Cref{cor:findexpander}.

\begin{proof}[ of \Cref{cor:findexpander}]
First, using a common folklore result, let $G_0$ be a bipartite subgraph $G$ with average degree at least $d(G)/2$. Then, we apply Lemma~\ref{lem:almostreg} to $G_0$ in order to find a (bipartite) $6$-almost-regular subgraph $G_1 \subseteq G_0$ with $d_1 \coloneqq d(G_1) \ge \frac{d(G)}{200 \log n}$. Now we apply \Cref{lem:findexpander} (with $G_1$ playing the role of $G$) to obtain a subgraph $G_2 \subseteq G_1$ which is an $(\eps, s)$-expander for some $s\geq \frac{d(G_2)}{(\log |V(G_2)|)^2}$, satisfying $d(G_2) \ge d_1 (1 - \frac{50}{\log \log |V(G_1)|})$ and $\delta(G_2) \ge \frac{d_1}{3}.$ \Cref{lem:findexpander} indeed applies since $d(G_1) \ge \frac{d(G)}{200 \log n} \ge (\log n)^2$. Moreover, since $|V(G_1)| \ge d(G_1) \ge (\log n)^2$, we have 
$d(G_2) \ge d_1 (1 - \frac{50}{\log \log |V(G_1)|}) \ge \frac{d_1}{2} \ge \frac{d(G)}{400 \log n}$.

Note that since $G_1$ is $6$-almost-regular, it has maximum degree at most $6 d_1$. Since $G_2$ is a subgraph of $G_1$, it also has maximum degree at most $6 d_1$, while, as noted above, its minimum degree is at least $\frac{d_1}{3}$, so $G_2$ is $18$-almost-regular. Thus, $G' \coloneqq G_2 \subseteq G$ is the desired $18$-almost-regular $(\eps,s)$-expander. This completes the proof of the lemma.
\end{proof}


\end{document}